\documentclass[10pt,english]{amsart}
\usepackage{amsfonts, amssymb, amsmath, amsthm, eucal, latexsym}

\usepackage[english]{babel}

\usepackage{pstricks}
\usepackage{graphicx}
\usepackage{psfrag}     
\usepackage{xy, xypic}
\theoremstyle{plain}
\newtheorem{theorem}{Theorem}[section]
\newtheorem{lemma}[theorem]{Lemma}
\newtheorem{coro}[theorem]{Corollary}
\newtheorem{prop}[theorem]{Proposition}
\newtheorem{conj}[theorem]{Conjecture}

\theoremstyle{definition}
\newtheorem{defi}[theorem]{Definition}

\theoremstyle{remark}
\newtheorem{rema}[theorem]{Remark}
\newtheorem{claim}[theorem]{Claim}

\theoremstyle{theorem}
\newtheorem{theo}[theorem]{Theorem}

\def\g{{\mathfrak{g}}}

\def\l{{\mathfrak{l}}}
\def\z{{\mathfrak{z}}}
\def\a{{\mathfrak{a}}}

\def\p{{\mathfrak{p}}}

\def\k{{\Bbbk}}

\def\ad{{\mathrm ad}\hskip .1em}   
\def\rk{{\rm rk}\hskip .1em}       
\def\rg{\ell}                      
\def\r{{\rm reg}}

\def\rs{{\rm reg,ss}}              
\def\ind{{\rm ind}\hskip .1em}     
\def\nil#1{\mathcal{N}(#1)}        
\def\orb#1{\mathcal{O}_{#1}}       
\def\slo{\mathcal{S}}              
\def\blw{B_e(\mathcal{S})}         

\def\til#1{\widetilde{#1}}     

\def\poi#1#2#3#4#5#6#7{\def\un{#5#6#7}\def\deux{#6#7}
\def\trois{#2#4} \def\cinq{#3#4#5}
\ifx\un\empty {#1}_{#2}{#3}{#1}_{#4} \else
\ifx\deux\empty {#5}(#1_{#2}){#3}{#5}(#1_{#4}) \else
\ifx\trois\empty {#5}_{#6}(#1){#3}{#5}_{#7}(#1) \else
{#5_{#6}}(#1_{#2}){#3}{#5_{#7}}(#1_{#4}) \fi \fi \fi}
\def\rond{\raisebox{.3mm}{\scriptsize$\circ$}}

\def\dv#1#2{\langle {#1},{#2}\rangle}

\def\tk#1#2{{#2}\otimes _{#1}}

\def\ec#1#2#3#4#5{\def\un{#3#4#5}\def\deux{#3#5}\def\trois{#3}
\def\four{#2#4#5}\def\five{#2#5}\def\six{#2}\def\seven{#3#4}
\def\eight{#2#4} \def\nine{#2#3#4}
\ifx\nine\empty {\rm #1}_{#5} \else
\ifx\un\empty {\rm #1}({\goth #2}) \else
\ifx\deux\empty {\rm #1}({\goth #2}_{#4}) \else
\ifx\trois\empty {\rm #1}_{#5}({\goth #2}_{#4}) \else
\ifx\four\empty {\rm #1}(#3) \else
\ifx\five\empty {\rm #1}(#3_{#4}) \else
\ifx\six\empty {\rm #1}_{#5}(#3_{#4}) \else
\ifx\seven\empty {\rm #1}_{#5} ({\goth#2})\else
\ifx\eight\empty {\rm #1}_{#5}({#3})
\fi \fi \fi \fi \fi \fi \fi \fi \fi}
\def\hec#1#2#3#4#5{\def\un{#3#4#5}\def\deux{#3#5}\def\trois{#3}
\def\four{#2#4#5}\def\five{#2#5}\def\six{#2}\def\seven{#3#4}
\def\eight{#2#4} \def\nine{#2#3#4}
\ifx\nine\empty \hat{{\rm #1}}_{#5} \else
\ifx\un\empty \hat{{\rm #1}}({\goth #2}) \else
\ifx\deux\empty \hat{{\rm #1}}({\goth #2}_{#4}) \else
\ifx\trois\empty \hat{{\rm #1}}_{#5}({\goth #2}_{#4}) \else
\ifx\four\empty \hat{{\rm #1}}(#3) \else
\ifx\five\empty \hat{{\rm #1}}(#3_{#4}) \else
\ifx\six\empty \hat{{\rm #1}}_{#5}(#3_{#4}) \else
\ifx\seven\empty \hat{{\rm #1}}_{#5} ({\goth#2})  \else
\ifx\eight\empty \hat{{\rm #1}}_{#5}({#3})
\fi \fi \fi \fi \fi \fi \fi \fi \fi}
\def\e#1#2{\ec {#1}#2{}{}{}}
\def\es#1#2{\ec {#1}{}{#2}{}{}}

\def\Bbb{\mathbb}
\def\goth{\mathfrak}
\def\cal{\mathcal}

\def\gi#1#2#3#4{\def\trois{#3#4} \def\quatre{#4}\def\cinq{#3}\ifx\trois\empty {\rm i}_{#1,{\goth #2}}
\else \ifx\quatre\empty {\rm i}_{#1_{#3},{\goth #2}} \else\ifx\cinq\empty {\rm i}_{#1,{\goth #2}_{#4}} \else {\rm i}_{#1_{#3},{\goth #2}_{#4}} \fi \fi \fi}
\def\j#1#2{\def\deux{#2} \ifx\deux\empty {\rm rk}\hskip .125em{{\goth #1}} \else {\rm rk}\hskip .125em{{\goth #1}_{#2}} \fi}
\def\aj#1#2{\def\deux{#2} \ifx\deux\empty {\rm j}_{{\goth #1}} \else {\rm j}_{{\goth #1}_{#2}} \fi}
\def\an#1#2{\def\deux{#2} \ifx\deux\empty {\cal O}_{#1} \else {\cal O}_{#1,#2} \fi }
\def\han#1#2{\def\deux{#2} \ifx\deux\empty {\hat{{\cal O}}}_{#1} \else {\hat{{\cal O}}}_{#1,#2} \fi }

\def\gg#1#2{\def\deux{#2} \ifx\deux\empty {\goth #1}\times {\goth #1}
\else {\goth #1}_{#2}\times {\goth #1}_{#2} \fi}
\def\ggr#1#2{\def\deux{#2} \ifx\deux\empty {\goth #1}_{\rs}\times {\goth #1}
\else ({\goth #1}_{#2})_{\rs}\times {\goth #1}_{#2} \fi}
\def\dim{{\rm dim}\hskip .125em}

\def\tr{{\rm tr}\hskip .125em}

\def\n{{\rm n}}
\def\s{{\rm s}}
\def\u{{\rm u}}

\def\b#1{{\mathrm {b}}_{{\mathfrak{#1}}}}

\title[The index of centralizers]{The index of centralizers of elements of reductive Lie algebras}

\author[J-Y Charbonnel]{Jean-Yves Charbonnel}

\address{Universit\'e Paris 7 - CNRS \\
Institut de Math\'ematiques de Jussieu \\
Th\'eorie des groupes \\
Case 7012 \\ B\^atiment Chevaleret \\
75205 Paris Cedex 13, France}

\email{jyc@math.jussieu.fr}

\author[A.~Moreau]{Anne Moreau}

\author[A. Moreau]{Anne Moreau}
\address{Anne Moreau, LMA\\
Boulevard Marie et Pierre Curie\\
86962 Futuroscope Chasseneuil Cedex, France}
\email{anne.moreau@math.univ-poitiers.fr}



\date\today

\begin{document}

\large

\maketitle

\begin{abstract}
For a finite dimensional complex Lie algebra, its index is the minimal dimension of
stabilizers for the coadjoint action. A famous conjecture due to A.G.~Elashvili says that
the index of the centralizer of an element of a reductive Lie algebra is equal to the
rank. That conjecture caught attention of several Lie theorists for years. It reduces to
the case of nilpotent elements. In~\cite{Pa1} and~\cite{Pa2}, D.I.~Panyushev proved the
conjecture for some classes of nilpotent elements (e.g.~regular, subregular and
spherical nilpotent elements). Then the conjecture has been proven for the classical Lie
algebras in~\cite{Y1} and checked with a computer programme for the exceptional
ones~\cite{De}. In this paper we give an almost general proof of that conjecture.
\end{abstract}

\section{Introduction}\label{int}
In this note $\k$ is an algebraically closed field of characteristic $0$.

\subsection{}\label{int1}
Let $\g$ be a finite dimensional Lie algebra over $\k$ and consider the coadjoint
representation of ${\goth g}$. By definition, the {\it index} of $\g$ is
the minimal dimension of stabilizers ${\goth g}^{x}$, $x\in {\goth g}^{*}$, for the
coadjoint representation:
$$ \ind {\goth g} \ := \ \min \{\dim \g^x; \ x \in \g^*\} \mbox{.}$$
The definition of the index goes back to Dixmier~\cite{Di1}.
It is a very important notion in representation theory and in invariant theory. By
Rosenlicht's theorem~\cite{Ro}, generic orbits of an arbitrary algebraic action of a
linear algebraic group on an irreducible algebraic variety are separated by rational
invariants; in particular, if $\g$ is an  algebraic Lie algebra,
$$\ind \g = \mathrm{deg\,tr}\,\k(\g^*)^{\g} ,$$
where $\k(\g^*)^{\g}$ is the field of $\g$-invariant rational functions on $\g^*$. The
index of a reductive algebra equals its rank. For an arbitrary Lie algebra, computing its
index seems to be a wild problem. However, there is a large number of interesting results
for several classes of nonreductive subalgebras of reductive Lie algebras. For instance,
parabolic subalgebras and their relatives as nilpotent radicals, seaweeds, are
considered in~\cite{Pa1}, \cite{TY1}, \cite{J}. The centralizers, or normalizers of
centralizers, of elements form another interesting class of such
subalgebras,~\cite{E1},~\cite{Pa1},~\cite{Mo2}. The last topic is closely related to
the theory of integrable Hamiltonian systems~\cite{Bol}. Let us precise this link.

From now on, $\g$ is supposed to be reductive. Denote by $G$ the adjoint group of
${\goth g}$. The symmetric algebra S$(\g)$ carries a natural Poisson structure.
By the so-called {\it argument shift method}, for $x$ in $\g ^{*}$,
we can construct a Poisson-commutative family $\mathcal{F}_{x}$ in
$\e Sg= \k[\g^*]$; see~\cite{MF} or Remark~\ref{rint3}.
It is generated by the derivatives of all orders in the direction $x\in\g^*$
of all elements of the algebra $\e Sg^{{\goth g}}$ of $\g$-invariants of S$(\g)$.
Moreover, if $G.x$ denotes the coadjoint orbit of $x\in \g^*$:

\begin{theo}[\cite{Bol}, Theorems 2.1 and 3.2]\label{tint1}
There is a Poisson-commutative family of polynomial functions on $\g^*$, constructed by
the argument shift method, such that its restriction to $G.x$ contains
$\frac{1}{2} \dim (G.x)$ algebraically independent functions if and only if
$\ind \g^x = \ind \g$.
\end{theo}

Denote by $\rk\g$ the rank of $\g$.
Motivated by the preceding result of Bolsinov, A.G.~Elashvili formulated a conjecture:
\begin{conj}[Elashvili]\label{cint}
Let  $\g$ be a reductive Lie algebra. Then $\ind\g^x=\rk\g$ for all $x\in\g^*$.
\end{conj}
Elashvili's conjecture also appears in the following problem:
Is the algebra S$(\g^x)^{\g^x}$ of invariants in $\es S{\g^{x}}$ under the adjoint action
a polynomial algebra? This question was formulated by A.~Premet
in~\cite[Conjecture 0.1]{PPY}. After that, O.~Yakimova discovered a
counterexample~\cite{Y3}, but the question remains very interesting.
As an example, under certain hypothesis and under the condition that Elashvili's
conjecture holds, the algebra of invariants S$(\g^x)^{\g^x}$ is polynomial in $\rk\g$
variables, \cite[Theorem 0.3]{PPY}.

During the last decade, Elashvili's conjecture caught attention of many invariant
theorists ~\cite{Pa1},~\cite{Ch},~\cite{Y1},~\cite{De}. To begin with, describe some easy
but useful reductions. Since the $\g$-modules $\g$ and $\g^*$ are isomorphic, it is
equivalent to prove Conjecture~\ref{cint} for centralizers of elements of $\g$. On the
other hand, by a result due to E.B.~Vinberg~\cite{Pa1}, the inequality
$\ind \g^x \geq \rk \g$ holds for all $x \in \g$. So it only remains to prove the opposite
one. Given $x\in \g$, let $x = x_{\s} + x_{\n}$ be its Jordan decomposition. Then
$\g^x = (\g^{x_{\s}})^{x_{\n}}$. The subalgebra $\g^{x_{\s}}$ is reductive of rank
$\rk\g$. Thus, the verification of Conjecture~\ref{cint} reduces to the case of
nilpotent elements. At last, one can clearly restrict oneself to the case of simple $\g$.

Review now the main results obtained so far on Elashvili's conjecture. If $x$ is regular,
then $\g^x$ is a commutative Lie algebra of dimension $\rk\g$. So, Conjecture~\ref{cint}
is obviously true in that case. Further, the conjecture is known for subregular nilpotent
elements and nilpotent elements of height 2 and 3,~\cite{Pa1},~\cite{Pa2}. Remind that
the {\it height} of a nilpotent element $e$ is the maximal integer $m$ such that
$(\ad e)^m \not= 0$. More recently, O.~Yakimova proved the conjecture in the classical
case~\cite{Y1}. To valid the conjecture in the exceptional types, W.~de~Graaf used the
computer programme \texttt{GAP}, see~\cite{De}. Since there are many nilpotent orbits in
the Lie algebras of exceptional type, it is difficult to present the results of such
computations in a concise way. In 2004, the first author published a case-free proof of
Conjecture~\ref{cint} applicable to all simple Lie algebras; see~\cite{Ch}.
Unfortunately, the argument in~\cite{Ch} has a gap in the final part of the proof which
was pointed out by L.~Rybnikov.

To summarize, so far, there is no conceptual proof of Conjecture~\ref{cint}.
Nevertheless, according to Yakimova's works and de~Graaf's works, we can claim:

\begin{theo}[\cite{Y1}, \cite{De}]\label{tint2}
Let  $\g$ be a reductive Lie algebra. Then $\ind\g^x=\rk\g$ for all $x\in\g^*$.
\end{theo}

Because of the importance of Elashvili's conjecture in invariant theory, it would be very
appreciated to find a general proof of Theorem~\ref{tint2} applicable to all
finite-dimensional simple Lie algebras. The proof we propose in this paper is fresh and
almost general. More precisely, it remains $7$ isolated cases; one nilpotent
orbit in type E$_7$ and six nilpotent orbits in type E$_8$ have to be considered
separately. For these 7 orbits, the use of \texttt{GAP} is unfortunately necessary.
In order to provide a complete proof of Theorem~\ref{tint2}, we include in this paper
the computations using \texttt{GAP} we made to deal with these remaining seven cases.

\subsection{Description of the paper}\label{int2}
Let us briefly  explain our approach. Denote by $\nil{\g}$ the nilpotent cone of $\g$.
As noticed previously, it suffices to prove $\ind \g^e=\rk \g$  for all $e$ in $\nil{\g}$.
If the equality holds for $e$, it does for all elements of $G.e$; we shortly say that
$G.e$ satisfies Elashvili's conjecture.

From a nilpotent orbit $\orb{\l}$ of a reductive factor $\l$ of a parabolic subalgebra of
$\g$, we can construct a nilpotent orbit of $\g$ having the same codimension in $\g$ as
$\orb{\l}$ in $\l$ and having other remarkable properties. The nilpotent orbits obtained
in such a way are called {\it induced}; the other ones are called {\it rigid}.
We refer the reader to Subsection~\ref{brn} for more precisions about this topic.
Using Bolsinov's criterion of Theorem~\ref{tint1}, we first prove Theorem~\ref{tint2}
for all induced nilpotent orbits and so the conjecture reduces to the case of rigid
nilpotent orbits. To deal with rigid nilpotent orbits, we use methods developed
in~\cite{Ch} by the first author, and resumed in~\cite{Mo1} by the second author, based
on nice properties of Slodowy slices of nilpotent orbits.\\

In more details, the paper is organized as follows:\\

We state in Section~\ref{pre} the necessary preliminary results.
In particular, we investigate in Subsection~\ref{cri} extensions of
Bolsinov's criterion and we establish an important result
(Theorem~\ref{tcri}) which will be used repeatedly in the sequel. We prove in
Section~\ref{cr} the conjecture for all induced nilpotent orbits (Theorem~\ref{tcr}) so
that Elashvili's conjecture reduces to the case of rigid nilpotent orbits
(Theorem~\ref{tcr}). From Section~\ref{ss}, we handle the rigid nilpotent orbits:
we introduce and study in Section~\ref{ss} a property (P) given by Definition~\ref{dss1}.
Then, in Section~\ref{ar}, we are able to deal with almost all rigid nilpotent orbits.
Still in Section~\ref{ar}, the remaining cases are dealt with set-apart by using a
different approach.

\subsection{Notations}\label{int3}
$\bullet$
If $E$ is a subset of a vector space $V$, we denote by span($E$) the vector subspace of
$V$ generated by $E$. The grassmanian of all $d$-dimensional subspaces of $V$ is denoted
by Gr$_d(V)$. By a {\it cone} of $V$, we mean a subset of $V$ invariant under the natural
action of $\k^{*}:=\k\setminus \{0\}$ and by a \emph{bicone} of $V\times V$ we mean a
subset of $V\times V$ invariant under the natural action of $\k^{*} \times \k^{*}$
on $V\times V$.

$\bullet$ From now on, we assume that $\g$ is semisimple of rank $\rg$
and we denote by $\dv ..$ the
Killing form of $\g$. We identify $\g$ to $\g^*$ through $\dv ..$. Unless otherwise
specified, the notion of orthogonality refers to the bilinear form $\dv ..$.

$\bullet$
Denote by $\e Sg^{\g}$ the algebra of $\g$-invariant elements of $\e Sg$.
Let $f_1,\ldots,f_{\rg}$ be homogeneous generators of $\e Sg^{\g}$ of degrees
$\poi d1{,\ldots,}{\rg}{}{}{}$ respectively. We choose the polynomials
$\poi f1{,\ldots,}{\rg}{}{}{}$ so that $\poi d1{\leq\cdots \leq }{\rg}{}{}{}$.
For $i=1,\ldots,\rg$ and $(x,y)\in\g \times \g$, we may consider a shift of $f_i$ in
direction $y$: $f_i(x+ty)$ where $t\in\k$. Expanding $f_i(x+ty)$ as a polynomial in $t$,
we obtain
\begin{eqnarray}\label{eq:fi}
f_i(x+ty)=\sum\limits_{m=0}^{d_i} f_{i}^{(m)} (x,y) t^m;  && \forall
(t,x,y)\in\k\times\g\times\g
\end{eqnarray}
where $y \mapsto (m!)f_{i}^{(m)}(x,y)$ is the differential at $x$ of $f_i$ of the order
$m$ in the direction $y$. The elements $f_{i}^{(m)}$ as defined by~(\ref{eq:fi}) are
invariant elements of $\tk {\k}{\e Sg}\e Sg$ under the diagonal action of $G$ on
$\gg g{}$. Note that $f_i^{(0)}(x,y)=f_i(x)$ while $f_i^{(d_i)}(x,y)=f_i(y)$  for all
$(x,y)\in \g\times \g$.

\begin{rema}\label{rint3}
The family
$\mathcal{F}_x  :=
\{f_{i}^{(m)}(x,.); \ 1 \leq i \leq \rg, 1 \leq m \leq d_i  \}$ for $x\in\g$,
is a Poisson-commutative family of $\e Sg$ by Mishchenko-Fomenko~\cite{MF}.
One says that the family $\mathcal{F}_x$ is constructed by
the \emph{argument shift method}.
\end{rema}

$\bullet$
Let $i \in\{1,\ldots,\rg\}$. For $x$ in $\g$, we denote by $\varphi_i(x)$ the element of
$\g$ satisfying $(\mathrm{d}f_i)_x(y)=f_{i}^{(1)}(x,y)=\dv {\varphi _{i}(x)}y$, for all
$y$ in ${\goth g}$.
Thereby, $\varphi_i$ is an invariant element of $\tk {\k}{\e Sg}\g$ under the canonical
action of $G$. We denote by $\varphi_{i}^{(m)}$, for
$0\leq m\leq d_i-1$, the elements of $\tk{\k}{\e Sg}\tk{\k}{\e Sg}\g$ defined by the
equality:
\begin{eqnarray}\label{eq:phi}
\varphi_i(x+ty) =\sum\limits_{ m=0}^{d_i-1} \varphi_i^{(m)}(x,y) t^m , &&
\forall (t,x,y)\in\k\times\g\times\g.
\end{eqnarray}

$\bullet$
For $x \in\g$, we denote by $\g^x= \ \{ y \in \g \ \vert \ [y,x]=0\}$ the centralizer of
$x$ in $\g$ and by $\z(\g^x)$ the center of $\g^x$. The set of regular elements of $\g$ is
$$\g_{\r} \ := \ \{ x\in \g \ \vert \ \dim \g^x=\rg \}$$
and we denote by ${\goth g}_{\rs}$ the set of regular semisimple elements of
${\goth g}$. Both $\g_{\r}$ and $\g_{\rs}$ are $G$-invariant dense open subsets of
${\goth g}$.

We denote by $C(x)$ the $G$-invariant cone generated by
$x$ and we denote by $x_{\s}$ and $x_{\n}$ the semisimple and nilpotent components of $x$
respectively.

$\bullet$
The nilpotent cone of $\g$ is $\nil{\g}$. As a rule, for $e \in \nil{\g}$, we choose an
$\mathfrak{sl}_2$-triple $(e,h,f)$ in $\g$ given by the Jacobson-Morozov theorem
\cite[Theorem~3.3.1]{CMa}. In particular, it satisfies the equalities:
$$ [h,e]=2e, \hspace{1cm} [e,f]=h, \hspace{1cm}  [h,f]=-2f $$
The action of $\ad h$ on $\g$ induces a $\mathbb{Z}$-grading:
$$ \g = \bigoplus_{i \in\mathbb{Z} } \g(i) \ , \
\g(i)=\{x\in \g \ \vert \ [h,x]=i x\}. $$
Recall that $e$, or  $G.e$, is said to be {\it even} if $\g(i)=0$ for odd $i$. Note
that $e\in \g(2)$, $f\in\g(-2)$ and that $\g^e$, $\z(\g^e)$ and $\g^f$ are all
$\ad h$-stable.

$\bullet$
All topological terms refer to the Zariski topology. If $Y$ is a subset of a topological
space $X$, we denote by $\overline{Y}$ the closure of $Y$ in $X$.

\subsection{Acknowledgments}
We would like to thank O.~Yakimova for her interest and useful discussions and
more particularly for bringing Bolsinov's paper to our attention.
We also thank A.G.~Elashvili for suggesting Lawther-Testerman's paper [LT08]
about the centers of centralizers of nilpotent elements.
We are grateful to Ilya Zakharevich for pointing out a mistake in Lemma~3.1 
in a previous version.

\setcounter{tocdepth}{1}
\tableofcontents

\section{Preliminary results}\label{pre}
We start in this section by reviewing some facts about the differentials of generators of
$\e Sg^{\g}$. Then, the goal of Subsection~\ref{cri} is Theorem~\ref{tcri}. We collect in
Subsection~\ref{brn} basic facts about induced nilpotent orbits.

\subsection{Differentials of generators of S$(\g)^{\g}$}\label{di}
According to subsection~\ref{int3}, the elements $\poi {\varphi }1{,\ldots,}{\rg}{}{}{}$
of $\tk {\k}{\e Sg}\g$ are the differentials of $\poi f1{,\ldots,}{\rg}{}{}{}$
respectively. Since $f_i(g(x))=f_i(x)$ for all $(x,g)\in\g \times  G$, the element
$\varphi _i(x)$ centralizes $x$ for all $x\in \g$. Moreover:

\begin{lemma}\label{ldi}
{\rm (i)\cite[Lemma 2.1]{Ri2}}
The elements $\poi x{}{,\ldots,}{}{\varphi }{1}{\rg}$ belong to $\z(\g^e)$.

{\rm (ii)\cite[Theorem 9]{Ko}}
The elements $\poi x{}{,\ldots,}{}{\varphi }{1}{\rg}$  are linearly independent elements
of $\g$ if and only if $x$ is regular. Moreover, if so,
$\poi x{}{,\ldots,}{}{\varphi }{1}{\rg}$  is a basis of $\g^x$.
\end{lemma}

We turn now to the elements $\varphi _{i}^{(m)}$, for $i=1,\ldots,\rg$ and
$0 \leq m \leq d_i-1$, defined in~Subsection~\ref{int3} by~(\ref{eq:phi}).
Recall that $d_i$ is the degree of the homogeneous polynomial $f_i$, for $i=1,\ldots,\rg$.
The integers $d_1-1,\ldots,d_{\rg}-1$ are thus the exponents of $\g$. By a classical
result~\cite[Ch.~V, \S5, Proposition 3]{Bou}, we have $\sum d_i=\b g$ where $\b g$ is the
dimension of Borel subalgebras of $\g$.
For $(x,y)$ in $\gg g{}$, we set:
\begin{eqnarray}\label{eq:V}
V_{x,y}:=\mathrm{span} \{\varphi_{i}^{(m)}(x,y) \ ; \  \ 1 \leq i \leq \rg,
0 \leq m \leq d_i-1  \} .
\end{eqnarray}
The subspaces $V_{x,y}$ will play a central role throughout the note.

\begin{rema}\label{rdi}
(1) For $(x,y)\in\g\times\g$, the dimension of $V_{x,y}$ is at most $\b g$ since
$\sum d_i= \b g$. Moreover, for all $(x,y)$ in a nonempty open subset
of $\g\times \g$, the equality holds~\cite{Bol}. Actually, in this note, we do not need
this observation.

(2) By Lemma~\ref{ldi}(ii), if $x$ is regular, then $\g^x$ is contained in $V_{x,y}$
for all $y \in\g$. In particular, if so, $\dim [x,V_{x,y}]=\dim V_{x,y}-\rg$.
\end{rema}

The subspaces $V_{x,y}$ were introduced and studied by Bolsinov in~\cite{Bol},
motivated by the maximality of Poisson-commutative families in $\e Sg$.
These subspaces have been recently exploited in~\cite{PY} and~\cite{CMo}.
The following results are mostly due to Bosinov,~\cite{Bol}. We refer to~\cite{PY} for a
more recent account about this topic. We present them in a slightly different way:

\begin{lemma}\label{l2di}
Let $(x,y)$ be in $\g _{\r}\times {\goth g}$.

{\rm (i)} The subspace $V_{x,y}$ of ${\goth g}$ is the sum of
the subspaces $\g^{x+ty}$ where $t$ runs through any nonempty open subset of $\k$ such
that $x+ty$ is regular for all $t$ in this subset.

{\rm (ii)} The subspace $\g^y+V_{x,y}$ is a totally isotropic subspace of ${\goth g}$
with respect to the Kirillov form $K_{y}$ on $\gg g{}$, $(v,w)\mapsto \dv y{[v,w]}$.
Furthermore, $\dim (\g^y+ V_{x,y})^{\perp} \geq \frac{1}{2}\dim G.y$.

{\rm (iii)} The subspaces $[x,V_{x,y}]$ and $[y,V_{x,y}]$ are equal.
\end{lemma}

\begin{proof}
(i) Let $O$ be a nonempty open subset of $\k$ such that $x+ty$ is regular for all $t$ in
$O$. Such an open subset does exist since $x$ is regular. Denote by $V_{O}$ the sum
of all the subspaces ${\goth g}^{x+ty}$ where $t$ runs through $O$. For all $t$ in $O$,
${\goth g}^{x+ty}$ is generated by $\poi {x+ty}{}{,\ldots,}{}{\varphi }{1}{\rg}$,
cf.~Lemma~\ref{ldi}(ii). As a consequence, $V_{O}$ is contained in $V_{x,y}$.
Conversely, for $i=1,\ldots,\rg$ and for $\poi t1{,\ldots,}{d_{i}}{}{}{}$ pairwise
different elements of $O$, $\varphi _{i}^{(m)}(x,y)$ is a linear combination of
$\varphi _{i}(x+t_{1}y),\ldots,\varphi _{i}(x+t_{d_{i}}y)$; hence
$\varphi _{i}^{(m)}(x,y)$ belongs to $V_{O}$. Thus $V_{x,y}$ is equal to $V_{O}$, whence
the assertion.

(ii) results from \cite[Proposition A4]{PY}. Notice that in (ii) the inequality is an
easy consequence of the first statement.

At last,~\cite[Lemma A2]{PY} gives us (iii).
\end{proof}

\noindent Let $\sigma $ and $\sigma_i$, for $i=1,\ldots,\rg$, be the maps
\begin{eqnarray*}
\begin{array}{ccl}
\g\times\g & \stackrel{\sigma}{\longrightarrow} &  \k^{\b g+\rg}\\
(x,y) & \longmapsto & (f_i^{(m)}(x,y))_{\hspace{-.2cm} 1\leq i\leq \rg,
\atop 0\leq m \leq d_i}
\end{array}, & \hspace{1cm}&
\begin{array}{ccl}
\g\times\g & \stackrel{\sigma_i}{\longrightarrow} &  \k^{d_i+1}\\
(x,y) & \longmapsto & (f_i^{(m)}(x,y))_{0\leq m \leq d_i}
\end{array}
\end{eqnarray*}
respectively, and denote by $\sigma '(x,y)$ and $\sigma _i'(x,y)$ the tangent map at
$(x,y)$ of $\sigma $ and $\sigma _i$ respectively. Then $\sigma _i'(x,y)$ is given by the
differentials of the $f_{i}^{(m)}$'s at $(x,y)$ and $\sigma '(x,y)$ is given by the
elements $\sigma _i'(x,y)$.

\begin{lemma}\label{l3di}
Let $(x,y)$ and $(v,w)$ be in $\gg g{}$.

{\rm (i)} For $i=1,\ldots,\rg$, $\sigma _i'(x,y)$ maps $(v,w)$ to
\begin{eqnarray*}
&& (\dv {\varphi_{i}(x)}v,\dv {\varphi_{i}^{(1)}(x,y)}v + \dv {\varphi_{i}^{(0)}(x,y)}w ,
\\  &&  \hspace{2.2cm} \ldots,
\dv {\varphi_{i}^{(d_i-1)}(x,y)}v +
\dv {\varphi_{i}^{(d_i-2)}(x,y)}w, \dv {\varphi_{i}(y)}w) .
\end{eqnarray*}

{\rm (ii)} Suppose that $\sigma'(x,y)(v,w)=0$. Then, for $w'$ in $\g$,
$\sigma'(x,y)(v,w')=0$ if and only if $w-w'$ is orthogonal to $V_{x,y}$.

{\rm (iii)} For $x \in\g_{\r}$, $\sigma '(x,y)(v,w')=0$ for some $w' \in \g$ if and
only if $v \in [x,{\goth g}]$.
\end{lemma}

\begin{proof}
(i) The verifications are easy and left to the reader.

(ii) Since $\sigma '(x,y)(v,w)=0$, $\sigma '(x,y)(v,w')=0$ if and only if
$\sigma '(x,y)(v,w-w')=0$ whence the statement by (i).

(iii) Suppose that $x$ is regular and suppose that $\sigma '(x,y)(v,w')=0$ for some
$w'\in\g$. Then by (i), $v$ is orthogonal to the elements
$\poi x{}{,\ldots,}{}{\varphi }{1}{\rg}$. So by Lemma~\ref{ldi}(ii), $v$ is orthogonal
to $\g ^{x}$. Since $\g^x$ is the orthogonal complement of $[x,\g]$ in $\g$, we deduce
that $v$ lies in $[x,\g]$. Conversely, since $\sigma (x,y)=\sigma (g(x),g(y))$ for
all $g$ in $G$, the element $([u,x],[u,y])$ belongs to the kernel of $\sigma '(x,y)$
for all $u\in\g$. So, the converse implication follows.
\end{proof}

\subsection{On Bolsinov's criterion}\label{cri}
Let $a$ be in $\g$ and denote by $\pi $ the map
\begin{eqnarray*}
{\goth g}  \times G.a & \stackrel{\pi}{\longrightarrow} & {\goth g} \times\k^{\b g+\rg} \\
(x,y) & \longmapsto & (x,\sigma (x,y)) .
\end{eqnarray*}

\begin{rema}\label{rcri}
Recall that the family $(\mathcal{F}_x)_{x\in\g}$ constructed by the argument shift method
consists of all elements $f_i^{(m)}(x,.)$ for $i=1,\ldots,\rg$ and $1 \leq m \leq d_i$,
see~Remark~\ref{rint3}. By definition of the morphism $\pi$, there is a family
constructed by the argument shift method whose restriction to $G.a$ contains
$\frac{1}{2}\dim G.a$ algebraically independent functions if and only if $\pi$ has a
fiber of dimension $\frac{1}{2}\dim G.a$.
\end{rema}

\noindent
In view of Theorem~\ref{tint1} and the above remark, we now concentrate on the fibers of
$\pi$. For $(x,y)\in  \g \times G.a$, denote by $F_{x,y}$ the fiber of $\pi $ at
$\pi (x,y)$:
$$F_{x,y} \ := \ \{x\} \times \{  y'\in G.a   \  \vert \ \sigma(x,y')=\sigma(x,y)\} .$$

\begin{lemma}\label{lcri}
Let $(x,y)$ be in $\g \times G.a $.

{\rm (i)} The irreducible components of $F_{x,y}$ have dimension at least
\sloppy \hbox{$\frac{1}{2}\dim G.a$}.

{\rm (ii)} The fiber $F_{x,y}$ has dimension $\frac{1}{2}\dim G.a$ if and only if any
irreducible component of $F_{x,y}$ contains an element $(x,y')$ such that
$(\g^{y'} + V_{x,y'})^{\perp}$ has dimension $\frac{1}{2}\dim G.a$.
\end{lemma}

\begin{proof}
We prove (i) and (ii) all together. The tangent space T$_{x,y'}(F_{x,y})$ of $F_{x,y}$
at $(x,y')$ in $F_{x,y}$ identifies to the subspace of elements $w$ of $[y',{\goth g}]$
such that $\sigma '(x,y')(0,w)=0$. Hence, by Lemma~\ref{l3di}(ii),
$$\mathrm{T}_{x,y'}(F_{x,y})= [y',{\goth g}] \cap V_{x,y'}^{\perp}
=(\g^{y'}+ V_{x,y'})^{\perp} ,$$
since $[y',\g]=(\g^{y'})^{\perp}$. But by Lemma~\ref{l2di}(ii),
$(\g^{y'}+ V_{x,y'})^{\perp}$ has dimension at least $\frac{1}{2}\dim G.a$; so does
$\mathrm{T}_{x,y'}(F_{x,y})$. This proves (i). Moreover, the equality holds if and
only if $(\g^{y'}+ V_{x,y'})^{\perp}$ has dimension $\frac{1}{2}\dim G.a$, whence the
statement (ii).
\end{proof}

\begin{theo}\label{tcri}
The following conditions are equivalent:
\begin{enumerate}
\item $\ind {\goth g}^{a}=\rg$;
\item $\pi $ has a fiber of dimension $\frac{1}{2}\dim G.a$;
\item there exists $(x,y)\in  {\goth g} \times G.a$ such that $(\g^y+ V_{x,y})^{\perp}$
has dimension $\frac{1}{2}\dim G.a$;
\item there exists $x$ in ${\goth g}_{\r}$ such that
$\dim ({\goth g}^{a}+V_{x,a})  = \frac{1}{2} (\dim {\goth g} + \dim {\goth g}^{a})$;
\item there exists $x$ in ${\goth g}_{\r}$ such that
$\dim V_{x,a} =  \frac{1}{2} \dim G.a + \rg$;
\item $\sigma ( \g \times \{a\})$ has dimension $\frac{1}{2}\dim G.a +\rg$.
\end{enumerate}
\end{theo}

\begin{proof}
By Theorem~\ref{tint1} and Remark~\ref{rcri}, we have (1)$\Leftrightarrow$(2).
Moreover, by Lemma~\ref{lcri}(ii), we have (2)$\Leftrightarrow$(3).

\noindent (3)$\Leftrightarrow$(4):
If (4) holds, so does (3).
Indeed, if so,
$$ \dim {\goth g} - \frac{1}{2}\dim G.a =
\frac{1}{2} (\dim {\goth g} + \dim {\goth g}^{a}) = \dim (\g^{a}+V_{x,a}). $$
Conversely, suppose that (3) holds. By Lemma~\ref{l2di}(ii),
$\g^y+ V_{x,y}$ has maximal dimension $\frac{1}{2} (\dim {\goth g}+\dim {\goth g}^{y})$.
So the same goes for all $(x,y)$ in a $G$-invariant nonempty open subset of
$\g \times G.a$. Hence, since the map $(x,y)\mapsto V_{x,y}$ is $G$-equivariant, there
exists $x$ in ${\goth g}_{\r}$ such that
$$\dim (V_{x,a}+{\goth g}^{a})  = \frac{1}{2} (\dim {\goth g} + \dim {\goth g}^{a}) .$$

\noindent (4)$\Leftrightarrow$(5):
Let $x$ be in ${\goth g}_{\r}$. By Lemma~\ref{l2di}(iii), $[x,V_{x,a}]=[a,V_{x,a}]$.
Hence ${\goth g}^{a}\cap V_{x,a}$ has dimension $\rg$ by Remark~\ref{rdi}(2).
As a consequence,
$$ \dim (\g^a+ V_{x,a}) = \dim {\goth g}^{a} + \dim V_{x,a} - \rg \mbox{,}$$
whence the equivalence.

\noindent (2)$\Leftrightarrow$(6): Suppose that (2) holds. By Lemma~\ref{lcri},
$\frac{1}{2}\dim G.a$ is the minimal dimension of the fibers of $\pi$.
So, $\pi (\g \times G.a )$ has dimension
$$ \dim {\goth g}+\dim G.a - \frac{1}{2}\dim G.a  =
\dim {\goth g}+\frac{1}{2}\dim G.a \mbox{.}$$
Denote by $\tau$ the restriction to $\pi (\g \times G.a)$ of the projection map
${\goth g}\times \k^{\b g+\rg} \to \k^{\b g+\rg}$. Then $\tau \rond \pi $ is the
restriction of $\sigma $ to $\g \times G.a$. Since $\sigma $ is a $G$-invariant map,
$\sigma ({\goth g}  \times \{a\}) = \sigma (\g \times G.a)$. Let
$(x,y) \in {\goth g}_{\rs} \times G.a$. The fiber of $\tau $ at $z=\sigma(x,y)$ is
$G.x$ since $x$ is a regular semisimple element of $\g$. Hence,
\begin{eqnarray*}\label{eq-26}
\dim \sigma ({\goth g} \times \{a\}) & = & \dim \pi (\g \times G.a) - (\dim {\goth g}-\rg)
                =  \frac{1}{2}\dim G.a+ \rg
\end{eqnarray*}
and we obtain (6).

Conversely, suppose that (6) holds. Then $\pi (\g \times G.a)$ has dimension
$\dim {\goth g}+\frac{1}{2}\dim G.a$ by the above equality. So the minimal dimension of
the fibers of $\pi $ is equal to
$$\dim {\goth g}+\dim G.a - (\dim {\goth g}+\frac{1}{2}\dim G.a )  =
        \frac{1}{2}\dim G.a  $$
and (2) holds.
\end{proof}

\subsection{Induced and rigid nilpotent orbits}\label{brn}
The definitions and results of this subsection are mostly extracted
from~\cite{Di1},~\cite{Di2},~\cite{LS} and~\cite{BoK}. We refer to~\cite{CMa} and
\cite{TY2} for recent surveys.

Let $\p$ be a proper parabolic subalgebra of $\g$ and let $\l$ be a reductive factor of
$\p$. We denote by $\p_{\u}$ the nilpotent radical of $\p$. Denote by $L$ the connected
closed subgroup of $G$ whose Lie algebra is $\ad \l$ and denote by $P$ the normalizer of
$\p$ in $G$.

\begin{theo}[\cite{CMa},Theorem 7.1.1]\label{tbrn}
Let $\mathcal{O}_{\l}$ be a nilpotent orbit of $\l$. There exists a unique nilpotent
orbit ${\cal O}_{\g}$ in ${\goth g}$ whose intersection with
$\mathcal{O}_{\l}+{\goth p}_{\u}$ is a dense open subset of
$\mathcal{O}_{\l}+{\goth p}_{\u}$. Moreover, the intersection of ${\cal O}_{\g}$ and
$\mathcal{O}_{\l}+{\goth p}_{\u}$ consists of a single $P$-orbit and
$\mathrm{codim}_{\g}({\cal O}_{\g})=\mathrm{codim} _{\l}({\cal O}_{\l})$.
\end{theo}

The orbit $\orb{\g}$ only depends on ${\goth l}$ and not on the choice of a parabolic
subalgebra ${\goth p}$ containing it~\cite[Theorem 7.1.3]{CMa}. By definition, the orbit
$\orb{\g}$ is called the \emph{induced orbit from ${\cal O}_{\l}$}; it is denoted by
Ind$_{\l}^{\g}(\orb{\l})$. If ${\cal O}_{\l}=0$, then we call ${\cal O}_{\g}$ a
\emph{Richardson orbit}. For example all even nilpotent orbits are
Richardson~\cite[Corollary 7.1.7]{CMa}. In turn, not all nilpotent orbits are induced
from another one. A nilpotent orbit which is not induced in a proper way from another one
is called {\it rigid}.

We shall say that $e \in \nil{\g}$ is an induced (respectively rigid) nilpotent element
of $\g$ if the $G$-orbit of $e$ is an induced (respectively rigid) nilpotent orbit of
$\g$. The following results are deeply linked to the properties of the sheets of
${\goth g}$ and the deformations of its $G$-orbits. We refer to~\cite{BoK} about these
notions.

\begin{theo}\label{t2brn}
{\rm (i)} Let $x$ be a non nilpotent element of ${\goth g}$ and let ${\cal O}_{\g}$ be the
induced nilpotent orbit from the adjoint orbit of $x_n$ in $\g^{x_s}$. Then
${\cal O}_{\g}$ is the unique nilpotent orbit contained in $\overline{C(x)}$ whose
dimension is $\dim G.x$. Furthermore,
$\overline{C(x)} \cap {\cal N}(\g)=\overline{{\cal O}_{\g}}$ and
$\overline{C(x)} \cap {\cal N}(\g)$ is the nullvariety in $\overline{C(x)}$
of $f_i$ where $i$ is an element of $\{1,\ldots,\ell\}$ such that $f_i(x)\not =0$.

{\rm (ii)} Conversely, if ${\cal O}_{\g}$ is an induced nilpotent orbit, there exists a
non nilpotent element $x$ of ${\goth g}$ such that
$\overline{C(x)} \cap {\cal N}(\g)=\overline{{\cal O}_{\g}}$.
\end{theo}

\begin{proof}
(i) Let ${\goth p}$ be a parabolic subalgebra of $\g$ having ${\goth g}^{x_{\s}}$ as a
Levi factor. Denote by ${\goth p}_{\u}$ its nilpotent radical and by $P$ the normalizer
of $\p$ in $G$. Let ${\cal O}'$ be the adjoint orbit of $x_{\n}$ in ${\goth g}^{x_{\s}}$.

\begin{claim} Let $C$ be the $P$-invariant closed cone generated by $x$ and
let $C_0$ be the subset of nilpotent elements of $C$.
Then $C=\k x_{\s}+\overline{{\cal O}'}+{\goth p}_{\u}$,
$C_{0}=\overline{{\cal O}'}+{\goth p}_{\u}$ and  $C_0$ is an irreducible subset of
dimension $\dim P(x)$.
\end{claim}

\begin{proof}
The subset $x_{\s}+\overline{{\cal O}'}+{\goth p}_{\u}$ is an irreducible closed subset
of ${\goth p}$ containing $P(x)$. Moreover, its dimension is equal to
$$ \dim {\cal O}'+\dim {\goth p}_{\u} = \dim {\goth g}^{x_{\s}}-\dim {\goth g}^{x}+
\dim {\goth p}_{\u} = \dim {\goth p}-\dim {\goth g}^{x} .$$
Since the closure of $P(x)$ and $x_{\s}+\overline{{\cal O}'}+{\goth p}_{\u}$ are both
irreducible subsets of $\g$, they coincide. As a consequence, the set
$\k x_{\s}+\overline{{\cal O}'}+{\goth p}_{\u}$ is contained in $C$.
Since the former set is clearly a closed conical subset of $\g$ containing $x$,
$C=\k x_{\s}+\overline{{\cal O}'}+{\goth p}_{\u}$. Then we deduce that
$C_{0}=\overline{{\cal O}'}+{\goth p}_{\u}$.
\end{proof}

Denote by $G\times _{P}{\goth g}$ the quotient of $G\times {\goth g}$ under the right
action of $P$ given by $(g,z).p:= (gp,p^{-1}(z))$. The map $(g,z)\mapsto g(z)$ from
$G\times {\goth g}$ to ${\goth g}$ factorizes through the quotient map from
$G\times {\goth g}$ to $G\times _{P}{\goth g}$. Since $G/P$ is a projective variety, the
so obtained map from $G\times _{P}{\goth g}$ to ${\goth g}$ is closed. Since $C$ and
$C_{0}$ are closed $P$-invariant subsets of ${\goth g}$, $G\times _{P}C$ and
$G\times _{P}C_{0}$ are closed subsets of $G\times _{P}{\goth g}$. Hence
$\overline{C(x)}=G(C)$ and $G(C_{0})$ is a closed subset of ${\goth g}$. So, by the
claim, the subset of nilpotent elements of $\overline{C(x)}$ is irreducible since
$C_{0}$ is irreducible. Since there are finitely many nilpotent orbits, the subset of
nilpotent elements of $\overline{C(x)}$ is the closure of one nilpotent orbit. Denote it
by $\tilde{{\cal O}}$ and prove $\tilde{{\cal O}}={\cal O}_{\g}$.

For all $k,l$ in $\{1,\ldots,\rg\}$, denote by $p_{k,l}$ the polynomial function
$$ p_{k,l} := f_{k}(x)^{d_{l}}f_{l}^{d_{k}}-f_{l}(x)^{d_{k}}f_{k}^{d_{l}} $$
Then $p_{k,l}$ is $G$-invariant and homogeneous of degree $d_{k}d_{l}$. Moreover
$p_{k,l}(x)=0$. As a consequence, $\overline{C(x)}$ is contained in the nullvariety of
the functions $p_{k,l}$, $1\leq k,l\leq \rg$. Hence the nullvariety of $f_{i}$ in
$\overline{C(x)}$ is contained in the nilpotent cone of ${\goth g}$ since it is the
nullvariety in ${\goth g}$ of the functions $\poi f1{,\ldots,}{\rg}{}{}{}$. Then
$\dim \tilde{{\cal O}} = \dim \overline{C(x)}-1 = \dim G.x$.
Since ${\cal O}'+{\goth p}_{\u}$ is contained in $\overline{C(x)}$,
Theorem~\ref{tbrn} tells us that ${\cal O}_{\g}$ is contained in $\overline{C(x)}$.
Moreover by Theorem~\ref{tbrn}, ${\cal O}_{\g}$ has dimension $\dim G.x$, whence
$\tilde{{\cal O}}={\cal O}_{\g}$. All statements of (i) are now clear.

(ii) By hypothesis,
${\cal O}_{\g}=\mathrm{Ind}_{\l}^{\g}(\orb{\l})$, where ${\goth l}$ is a proper Levi
subalgebra of ${\goth g}$ and ${\cal O}_{\l}$ a nilpotent orbit in ${\goth l}$. Let
$x_{\s}$ be an element of the center of ${\goth l}$ such that
${\goth g}^{x_{\s}}={\goth l}$, let $x_{\n}$ be an element of ${\cal O}_{\l}$ and set
$x=x_{\s}+x_{\n}$. Since $\l$ is a proper subalgebra, the element $x$ is not nilpotent.
So by (i), the subset of nilpotent elements of $\overline{C(x)}$ is the closure of
${\cal O}_{\g}$.
\end{proof}

\section{Proof of Theorem~\ref{tint2} for induced nilpotent orbits}\label{cr}
Let $e$ be an induced nilpotent element. Let $x$ be a non nilpotent element of
${\goth g}$ such that $\overline{C(x)} \cap {\cal N}(\g)= \overline{G.e}$. Such an
element does exist by Theorem~\ref{t2brn}(ii).

\begin{theo}\label{tcr}
Assume that $\ind {\goth a}^{x}=\rk {\goth a}$ for all reductive subalgebras ${\goth a}$
strictly contained in $\g$ and for all $x$ in ${\goth a}$. Then for all induced
nilpotent orbits $\orb {\g}$ in ${\goth g}$ and for all $e$ in $\orb {\g}$,
$\ind\g^e=\rg$.
\end{theo}

\begin{proof}
Let ${\cal O}_{\g}$ be an induced nilpotent orbit and let $e$ be in ${\cal O}_{\g}$.
Using Theorem 2.9(ii), we let $x$ be a non nilpotent element of ${\goth g}$ such
that $\overline{C(x)} \cap {\cal N}(\g)= \overline{{\cal O}_{\g}}$.
Since $x$ is not nilpotent, ${\goth g}^{x}$ is the
centralizer in the reductive Lie algebra ${\goth g}^{x_{\s}}$ of the nilpotent element
$x_{\n}$ of ${\goth g}^{x_{\s}}$. Since ${\goth g}^{x_{\s}}$ is strictly contained in
${\goth g}$ and has rank $\rg$, the index of ${\goth g}^{x}$ is equal to $\rg$ by
hypothesis. Besides, by Theorem 2.7, (1)$\Rightarrow$(6), applied to $x$,
$$ \dim \sigma ({\goth g}\times \{x\}) = \frac{1}{2}\dim G.x + \ell.$$
Since $\sigma $ is $G$-invariant,
$\sigma ({\goth g}\times \{x\})=\sigma ({\goth g}\times G.x)$.
Hence for all $z$ in a dense subset of $\sigma ({\goth g}\times G.x)$, the fiber of the
restriction of $\sigma $ to ${\goth g}\times G.x$ at $z$ has minimal dimension
\begin{eqnarray*}
\dim {\goth g}+ \dim G.x - (\frac{1}{2}\dim G.x + \ell) =
\dim {\goth g}+\frac{1}{2}\dim G.x - \ell .
\end{eqnarray*}
Denote by $Z$ the closure of $\sigma ({\goth g}\times \overline{C(x)})$ in
$\k ^{\underline{d}}$.
We deduce from the above equality that $Z$ has dimension
\begin{eqnarray*}
    \dim {\goth g}+\dim C(x) - (\dim {\goth g}+\frac{1}{2}\dim G.x - \ell) & = &
    \dim C(x) - \frac{1}{2} \dim G.x + \ell \\
    & = & \frac{1}{2} \dim G.e + \ell +1,
\end{eqnarray*}
since $\dim C(x)=\dim G.x+1=\dim G.e+1$.

Let $i$ be in $\{1,\ldots,\rg\}$ such that $f_i(x)\not=0$. For $z\in \k^{\underline{d}}$,
we write
$z=(z_{i,j})_{\hspace{-.15cm} 1\leq i\leq  \rg \atop 0\leq j \leq d_i}$ its
coordinates. Let $\mathcal{V}_{i}$ be the nullvariety in 
$\sigma ({\goth g}\times \overline{C(x)})$ of the coordinate $z_{i,d_{i}}$. Then 
$\mathcal{V}_{i}$ is not empty. Since $\sigma ({\goth g}\times \overline{C(x)})$ is an
irreducible constructible subset of $\k^{\underline{d}}$ and since $z_{i,d_{i}}$ is not 
identically zero on $\sigma ({\goth g}\times \overline{C(x)})$, $\mathcal{V}_{i}$ has
dimension $\frac{1}{2}\dim G.e+\rg$. By Theorem 2.9(i),
the nullvariety of $f_{i}$ in $\overline{C(x)}$ is equal to $\overline{G.e}$. Hence
$$ {\goth g}\times \overline{G.e} = 
\sigma ^{-1}(\mathcal{V}_{i})\cap ({\goth g}\times \overline{C(x)}) $$
So $\sigma ({\goth g}\times G.e)$ is equal to ${\cal V}_{i}$ and has dimension 
$\frac{1}{2}\dim G.e+\rg$. Then by Theorem 2.7,
(6)$\Rightarrow$(1), the index of ${\goth g}^{e}$ is equal to $\rg$.
\end{proof}

From that point, our goal is to prove Theorem~\ref{tint2} for rigid nilpotent elements;
Theorem~\ref{tcr} tells us that this is enough to complete the proof.

\section{The Slodowy slice and the property (P)}\label{ss}
In this section, we introduce a property (P) in Definition~\ref{dss1} and we prove that
$e \in \nil{\g}$ has Property (P) if and only if $\ind\g^e=\rg$ (Theorem~\ref{tss2}).
Then, we will show in the next section that all rigid nilpotent orbits of $\g$ but seven
orbits (one in the type ${\mathrm {E}}_{7}$ and six in the type ${\mathrm {E}}_{8}$)
do have Property (P).

\subsection{Blowing up of $\slo$}\label{ss1}
Let $e$ be a nilpotent element of $\g$ and consider an $\mathfrak{sl}_2$-triple
$(e,h,f)$ containing $e$ as in Subsection~\ref{int3}.
The {\it Slodowy slice} is the
affine subspace $\slo:=e+\g^f$ of $\g$ which is a transverse variety to the adjoint orbit
$G.e$. Denote by $\blw$ the blowing up of $\slo$ centered at $e$ and let
$p : \blw \to \slo$ be the canonical morphism. The variety $\slo$ is smooth and
$p^{-1}(e)$ is a smooth irreducible hypersurface of $\blw$. The use of the blowing-up
$\blw$ for the computation of the index was initiated by the first author in~\cite{Ch}
and resumed by the second author in~\cite{Mo1}. Here, we use again this technique to
study the index of $\g^e$. Describe first the main tools extracted from~\cite{Ch} we need.

For $Y$ an open subset of $\blw$, we denote by $\k[Y]$ the algebra of regular functions
on $Y$. By \cite[Th\'eor\`eme 3.3]{Ch}, we have:

\begin{theo}\label{tss1}
The following two  assertions are equivalent:

{\rm (A)} the equality  $\ind\g^e=\rg$ holds,

{\rm (B)} there exists an affine open subset $Y \subset \blw$ such that
$Y\cap p^{-1}(e)\not=\emptyset$ and satisfying the following property:
\begin{quotation}
for any regular map $\varphi \in\k[Y]\otimes_{\k} \g$ such that
$\varphi(x)\in [{\goth g},p(x)]$ for all $x\in Y$, there exists
$\psi \in\k[Y]\otimes_{\k} \g$ such that $\varphi (x)=[\psi (x),p (x)]$ for all $x\in Y$.
\end{quotation}
\end{theo}

An open subset $\Omega \subset \blw$ is called {\it a big open subset} if
$\blw\setminus \Omega$ has codimension  at least $2$ in $\blw$. As explained in
\cite[Section 2]{Ch}, there exists a big open subset $\Omega $ of $\blw$ and a regular map
$$\alpha :\Omega \to \mathrm{Gr}_{\rg}(\g) $$
such that $\alpha (x)=\g^{p(x)}$ if $p(x)$ is regular. Furthermore, the map $\alpha $ is
uniquely defined by this condition. In fact, this result is a consequence of
\cite[Ch.~VI, Theorem 1]{Sh}.
From now on, $\alpha$ stands for the so-defined map.
Since $p^{-1}(e)$ is an hypersurface and since $\Omega$ is a big open subset of $\blw$,
note that $\Omega \cap p^{-1}(e) $ is a nonempty set. In addition,
$\alpha(x)\subset \g^{p(x)}$ for all $x\in \Omega $.

\begin{defi}\label{dss1}
We say that $e$ \emph{has Property} (P) if $\z(\g^e)\subset \alpha (x)$ for all $x$ in
$\Omega \cap  p^{-1}(e) $.
\end{defi}

\begin{rema}\label{rss1}
Suppose that $e$ is regular. Then $\g^e$ is a commutative algebra, i.e.~$\z(\g^e)=\g^e$.
If $x \in \Omega  \cap p^{-1}(e)$, then $\alpha (x)=\g^{e}$ since $p(x)=e$ is regular in
this case. On the other hand, $\ind\g^e=\dim \g^e=\rg$ since $e$ is regular. So $e$ has
Property (P) and $\ind\g^e=\rg$.
\end{rema}

\subsection{On the property (P)}\label{ss2}
This subsection aims to show: Property (P) holds for $e$ if and only if
$\ind\g^e=\rg$. As a consequence of Remark~\ref{rss1}, we can (and will) assume that $e$
is a nonregular nilpotent element of $\g$. As a first step, we will state in
Corollary~\ref{c2ss2} that, if (P) holds, then so does the assertion (B) of
Theorem~\ref{tss1}.

Let $L_{\g}$ be the $\e Sg$-submodule of $\varphi \in \tk {\k}{\e Sg}\g$ satisfying
$[\varphi (x),x]=0$ for all $x$ in $\g$. It is known that $L_{{\goth g}}$ is a free
module of basis $\varphi_1,\ldots,\varphi_{\rg}$, cf.~\cite{Di3}. We investigate an
analogous property for the Slodowy slice $\slo=e+\g^f$. We denote by $\slo_{\r}$ the
intersection of $\slo$ and $\g_{\r}$. As $e$ is nonregular, the set
$(\slo\setminus \slo_{\r})$ contains $e$.

\begin{lemma}\label{lss2}
The set $\slo\setminus \slo_{\r}$ has codimension $3$ in $\slo$ and each irreducible
component of $\slo\setminus \slo_{\r}$ contains $e$.
\end{lemma}

\begin{proof}
Let us consider the morphism
\begin{eqnarray*}
G\times \slo & \longrightarrow & {\goth g} \\
(g,x) & \longmapsto & g(x)
\end{eqnarray*}
By a Slodowy's result~\cite{Sl}, this morphism is a smooth morphism. So its fibers are
equidimensional of dimension $\dim {\goth g}^{f}$. In addition, by~\cite{V},
$\g\setminus \g_{\r}$ is a $G$-invariant equidimensional closed subset of ${\goth g}$ of
codimension 3. Hence $\slo\setminus {\slo}_{\r}$ is an equidimensional closed subset of
$\slo$ of codimension $3$.

Denoting by $t\mapsto g(t)$ the one parameter subgroup of $G$ generated by $\ad h$,
$\slo$ and $\slo \setminus \slo _{\r}$ are stable under the action of $t^{-2}g(t)$
for all $t$ in $\k ^{*}$. Furthermore, for all $x$ in $\slo$, $t^{-2}g(t)(x)$ goes to
$e$ when $t$ goes to $\infty $, whence the lemma.
\end{proof}

Denote by $\k[\slo]$ the algebra of regular functions on $\slo$ and denote by $L_{\slo}$
the $\k[\slo]$-submodule of $\varphi \in  \tk {\k}{\k[\slo]}\g$ satisfying
$[\varphi (x),x]=0$ for all $x$ in $\slo$.

\begin{lemma}\label{l2ss2}
The module $L_{\slo}$ is a free module of basis
$\varphi _{1}  \vert _{\slo},\ldots,\varphi _{\rg }\vert _{\slo} $ where
$\varphi _{i} \vert _{\slo}$ is the restriction to $\slo$ of $\varphi _{i}$
for $i=1,\ldots,\rg$.
\end{lemma}

\begin{proof}
Let $\varphi $ be in $L_{\slo}$. There are regular functions $a_{1},\ldots,a_{\rg}$ on
$\slo_{\r}$ satisfying
$$\varphi (x) = a_{1}(x)\varphi _{1} \vert _{\slo} (x) + \cdots +
a_{\rg}(x)\varphi _{\rg} \vert _{\slo} (x) $$
for all $x \in {\slo}_{\r}$, by Lemma~\ref{ldi}(ii). By Lemma~\ref{lss2},
$\slo\setminus\slo_{\r}$ has codimension $3$ in $\slo$. Hence $a_{1},\ldots,a_{\rg}$ have
polynomial extensions to $\slo$ since $\slo$ is normal. So the maps
$\varphi_{1} \vert _{\slo},\ldots,\varphi_{\rg} \vert _{\slo} $
generate $L_{\slo}$.
Moreover, by Lemma~\ref{ldi}(ii) for all $x \in {\slo}_{\r}$,
$\varphi _{1}(x),\ldots,\varphi _{\rg}(x)$ are linearly independent, whence the statement.
\end{proof}

The following proposition accounts for an important step to interpret Assertion (B) of
Theorem~\ref{tss1}:

\begin{prop}\label{pss2}
Let $\varphi $ be in $\tk {\k}{\k[\slo]}\g$ such that $\varphi(x)\in [{\goth g},x]$
for all $x$ in a nonempty open subset of $\g$. Then there exists a polynomial map
$\psi\in \tk {\k}{\k[\slo]}\g$ such that $\varphi(x)=[\psi(x),x]$ for all $x \in \slo$.
\end{prop}

\begin{proof}
Since ${\goth g}^{x}$ is the orthogonal complement of $[x,{\goth g}]$ in ${\goth g}$,
our hypothesis says that $\varphi(x)$ is orthogonal to $\g^x$ for all $x$ in a nonempty
open subset $\slo'$ of $\slo $. The intersection $\slo '\cap \slo_{\r}$ is not empty;
so by Lemma~\ref{ldi}(ii), $\dv {\varphi(x)}{\varphi _{i} \vert _{\slo} (x)}=0$ for all
$i=1,\ldots,\rg$ and for all $x\in \slo'\cap \slo_{\r}$. Therefore, by continuity,
$\dv {\varphi(x)}{\varphi _{i} \vert _{\slo} (x)}=0$ for all $i=1,\ldots,\rg$ and
all $x\in\slo$. Hence $\varphi(x) \in [x,{\goth g}]$ for all $x\in {\slo}_{\r}$
by Lemma~\ref{ldi}(ii) again. Consequently by Lemma \ref{lss2}, Lemma \ref{l2ss2} and
the proof of the main theorem of~\cite{Di3}, there exists an element $\psi \in
\tk {\k}{\k[\slo]}{\goth g}$ which satisfies the condition of the proposition.
\end{proof}

Let $\poi u1{,\ldots,}{m}{}{}{}$ be a basis of ${\goth g}^{f}$ and let
$u_{1}^{*},\ldots,u_{m}^{*}$ be the corresponding coordinate system of $\slo=e+\g^f$.
There is an affine open subset $Y \subset \blw$  with $Y \cap p ^{-1}(e)\not=\emptyset$
such that $\k[Y]$ is the set of linear combinations of monomials in
$(u_{1}^{*})^{-1},u_{1}^{*},\ldots,u_{m}^{*}$ whose total degree is nonnegative. In
particular, we have a global coordinates system $u_{1}^{*},v_{2}^{*},\ldots,v_{m}^{*}$ on
$Y$ satisfying the relations:
\begin{eqnarray}\label{eq:vi}
u_{2}^{*}= u_{1}^{*}v_{2}^{*} & ,  \ldots , &  u_{m}^{*} = u_{1}^{*}v_{m}^{*} .
\end{eqnarray}
Note that, for $x \in Y$, we so have:
$p(x)=e+u_{1}^{*}(x) (u_1 +  v_{2}^{*}(x) u_2 + \cdots + v_{m}^{*}(x) u_m)$.
So, the image of $Y$ by $p$
is the union of $\{e\}$ and the complementary in $\slo$ of the nullvariety of $u_{1}^{*}$.
Let $Y'$ be an affine open subset of $Y$ contained in $\Omega $ and having a nonempty
intersection with $p ^{-1}(e)$. Denote by $L_{Y'}$ the set of regular maps
$\varphi $ from $Y'$ to ${\goth g}$ satisfying $[\varphi (x),p(x)]=0$ for all $x\in Y'$.

\begin{lemma}\label{l3ss2}
Suppose that $e$ has Property $(\mathrm{P})$. For each $z\in\z({\g}^e)$, there exists
$\psi _{z}\in\tk {\k}{\k[Y']}{\goth g}$ such that $z - u_{1}^{*} \psi _{z}$ belongs to
$L_{Y'}$.
\end{lemma}

\begin{proof}
Let $z$ be in $\z(\g^e)$. Since $Y'\subset\Omega $, for each $y\in Y'$, there exists an
affine open subset $U_{y}$ of $Y'$ containing $y$ and regular maps
$\poi {\nu }1{,\ldots,}{\rg}{}{}{}$ from $U_{y}$ to ${\goth g}$ such that
$\poi x{}{,\ldots,}{}{\nu }{1}{\rg}$ is a basis of $\alpha (x)$ for all $x\in U_{y}$.
Let $y$ be in $Y'$. We consider two cases:

(1) Suppose $p(y)=e$.\\
Since $e$ has Property (P), there exist regular functions
$\poi a1{,\ldots,}{\rg}{}{}{}$ on $U_{y}$ satisfying
$$ z = a_{1}(x)\nu _{1}(x) + \cdots + a_{\rg}(x)\nu _{\rg}(x) \mbox{,}$$
for all $x \in U_{y} \cap p^{-1}(e)$. The intersection $U_{y} \cap p ^{-1}(e)$ is the set
of zeroes of $u_{1}^{*}$ in $U_{y}$. So there exists a regular map $\psi $ from $U_{y}$ to
${\goth g}$ which satisfies the equality:
$$z - u_{1}^{*}\psi = a_{1}\nu _{1} + \cdots + a_{\rg}\nu _{\rg}.$$
Hence $[z-u_{1}^{*}(x)\psi(x),p(x)]=0$ for all $x\in U_y$ since $\alpha(x)$ is contained
in $\g^{p(x)}$ for all $x\in\Omega$.

(2) Suppose $p(y)\not=e$.\\
Then we can assume that $U_{y}\cap p ^{-1}(e)=\emptyset$ and the map
$\psi=(u_{1}^{*})^{-1}z$ satisfies the condition:
$[z-u_{1}^{*}(x)\psi (x),p(x)]=0$ for all $x \in U_y$.

In  both cases (1) or (2), we have found a regular map $\psi _y$ from $U_y$
to $\g$ satisfying: $[z-(u_{1}^{*} \psi_y)(x),p(x)]=0$ for all $x \in U_y$.

Let $\poi y1{,\ldots,}{k}{}{}{}$ be in $Y'$ such that the open subsets
$\poi U{y_{1}}{,\ldots,}{y_{k}}{}{}{}$ cover $Y'$. For $i=1,\ldots,k$, we denote by
$\psi _{i}$ a regular map from $U_{y_{i}}$ to ${\goth g}$ such that
$z-u_{1}^{*}\psi _{i}$ is in $\Gamma (U_{y_{i}},{\cal L})$ where ${\cal L}$ is the
localization of $L_{Y'}$ on $Y'$. Then for $i,j=1,\ldots,k$, $\psi _{i}-\psi _{j}$ is in
$\Gamma (U_{y_{i}}\cap U_{y_{j}},{\cal L})$. Since $Y'$ is affine, H$^{1}(Y',{\cal L})=0$.
So, for $i=1,\ldots,l$, there exists $\til{\psi} _{i}$ in
$\Gamma (U_{y_{i}},{\cal L})$ such that $\til{\psi} _{i}-\til{\psi} _{j}$ is equal to
$\psi _{i}-\psi _{j}$ on $U_{y_{i}}\cap U_{y_{j}}$ for all $i,j$.
Then there exists a well-defined map $\psi _z$ from $Y'$ to ${\goth g}$ whose restriction
to $U_{y_{i}}$ is equal to $\psi _{i}-\til{\psi} _{i}$ for all $i$, and such that
$z-u_{1}^{*}\psi _z $ belongs to $L_{Y'}$. Finally, the map $\psi _{z}$ verifies the
required property.
\end{proof}

Let $z$ be in $\z({\g}^e)$. We denote by $\varphi_{z}$ the regular map from $Y$ to $\g$
defined by:
\begin{eqnarray}\label{eq:phi-z}
\varphi_{z}(x) =  [z,u_{1}] + v_{2}^{*}(x) [z,u_{2}] + \cdots +
    v_{m}^{*}(x) [z,u_{m}] , \ \textrm{ for all } x \in Y .
\end{eqnarray}

\begin{coro}\label{css2}
Suppose that $e$ has Property $(\mathrm{P})$ and let $z$ be in $\z({\g}^e)$.
There exists $\psi _{z}$ in $\tk {\k}{\k[Y']}{\goth g}$ such that
$\varphi _{z}(x)=[\psi _{z}(x),p(x)]$ for all $x \in Y'$.
\end{coro}

\begin{proof}
By Lemma \ref{l3ss2}, there exists $\psi _{z}$ in $\tk {\k}{\k[Y']}{\goth g}$ such that
$z-u_{1}^{*}\psi_{z}$ is in $L_{Y'}$. Then
$$u_{1}^{*} \varphi_{z}(x) = [z,p (x)]  = [z-u_{1}^{*} \psi_{z}(x),p (x)] +
u_{1}^{*} [\psi_{z}(x),p (x)] \mbox{,}$$
for all $x\in Y'$. So the map $\psi_{z}$
is convenient, since $u_{1}^{*}$ is not identically zero on $Y'$.
\end{proof}

The following lemma is easy but helpful for Proposition~\ref{p2ss2}:

\begin{lemma}\label{l4ss2}
Let $v$ be in ${\goth g}^{e}$. Then, $v$ belongs to $\z({\g}^e)$ if and only
if $[v,{\g}^{f}] \subset [e,{\goth g}]$.
\end{lemma}

\begin{proof}
Since $[x,{\goth g}]$ is the orthogonal complement of ${\goth g}^{x}$ in $\g$ for all
$x\in\g$, we have:
\begin{eqnarray*}
[v,{\g}^{f}] \subset [e,{\goth g}] \iff \dv {[v,\g^f]}{\g^e}=0
\iff \dv {[v,\g^e]}{\g^f}=0 \iff [v,\g^e] \subset [f,\g].
\end{eqnarray*}
But $\g$ is the direct sum of $\g^e$ and $[f,\g]$ and $[v,\g^e]$ is contained in $\g^e$
since $v\in\g^e$. Hence $[v,{\g}^{f}]$ is contained in $[e,{\goth g}]$ if and only if
$v$ is in $\z(\g^e)$.
\end{proof}

\begin{prop}\label{p2ss2}
Suppose that $e$ has Property $(\mathrm{P})$ and let $\varphi $ be in
$\tk {\k}{\k[Y]}{\goth g}$ such that \hbox{$\varphi(x) \in [{\goth g},p (x)]$} for all
$x \in Y$. Then there exists $\psi $ in $\tk {\k}{\k[Y']}{\goth g}$ such that
$\varphi(x)=[\psi(x),p(x)]$ for all $x \in Y'$.
\end{prop}

\begin{proof}
Since $\varphi$ is a regular map from $Y$ to ${\goth g}$, there is a nonnegative
integer $d$ and $\til{\varphi}\in \tk {\k}{\k[\slo]}\g$ such that
\begin{eqnarray}\label{eq:phi-1}
(u_{1}^{*})^{d}(x) \varphi (x) = (\til{\varphi} \rond p )(x), \ \forall x \in Y
\end{eqnarray}
and $\til{\varphi}$ is a linear combination of monomials in $u_{1}^{*},\ldots,u_{m}^{*}$
whose total degree is at least $d$.
By hypothesis on $\varphi$, we deduce that for all $x\in \slo$ such that
$u_{1}^{*}(x)\not=0$, $\til{\varphi}(x)$ is in $[{\goth g},x]$. Hence by Proposition
\ref{pss2}, there exists $\til{\psi}$ in $\tk {k}{\k[\slo]}\g$ such that
$\til{\varphi}(x)=[\til{\psi}(x),x]$ for all $x \in \slo$.

Denote by $\til{\psi}'$ the sum of monomials of degree at least $d$ in $\til{\psi}$
and denote by $\psi'$ the element of $\tk {\k}{\k[Y]}\g$ satisfying
\begin{eqnarray}\label{eq:psi-1}
(u_{1}^{*})^{d}(x)\psi'(x)=(\til{\psi}' \rond p)(x), \ \forall x\in Y.
\end{eqnarray}
Then we set, for $x \in Y$, $\varphi'(x):=\varphi(x)-[\psi'(x),p (x)]$.
We have to prove the existence of an element  $\psi''$ in $\tk {\k}{\k[Y']}{\goth g}$
such that $\varphi'(x)=[\psi''(x),p (x)]$ for all $x\in Y'$.

$\bullet$
If $d=0$, then $\varphi=\til{\varphi}\rond p$, $\psi'=\psi$ and $\varphi'=0$; so
$\psi'$ is convenient in that case.

$\bullet$
If $d=1$, we can write
\begin{eqnarray*}
u_{1}^{*}(x)\varphi(x)=\til{\varphi}(p(x))
        =[\til{\psi}(p(x)),e+u_{1}^{*}(x)(u_1+v_{2}^{*}(x)u_2+\cdots+v_{m}^{*}(x)u_m)],
\end{eqnarray*}
for all $x\in Y$, whence we deduce
\begin{eqnarray*}
u_{1}^{*}(x)(\varphi (x)-[\psi '(x),p(x)]) =
[\tilde{\psi }(e),e+u_{1}^{*}(x)(u_1+v_{2}^{*}(x)u_2+\cdots+v_{m}^{*}(x)u_m)]
\end{eqnarray*}
for all $x\in Y$. Hence $\tilde{\psi }(e)$ belongs to ${\goth g}^{e}$ and
$[\til{\psi}(e),u_{i}] \in [e,{\goth g}]$ for all $i=1,\ldots,m$,
since $\varphi(x)\in [e,{\goth g}]$ for all $x\in Y \cap p ^{-1}(e)$. Then
$\til{\psi}(e)$ is in $\z({\goth g}^{e})$ by Lemma~\ref{l4ss2}. So by
Corollary~\ref{css2}, $\varphi'$ has the desired property.

$\bullet$
Suppose $d>1$. For $\underline{i}=(i_{1},\ldots,i_{m})\in{\Bbb N}^{m}$, we set
$\vert \underline{i} \vert:= i_{1}+\cdots+i_{m}$ and we denote by
$\psi _{\underline{i}}$ the coefficient of
$(u_{1}^{*})^{i_{1}}\cdots (u_{m}^{*})^{i_{m}}$ in $\til{\psi}$. By Corollary \ref{css2},
it suffices to prove:
$$\left\{\begin{array}{ll}
\psi_{\underline{i}}=0 & \textrm{ if } \vert \underline{i} \vert < d-1 ; \\
\psi_{\underline{i}}\in\z({\g}^e) & \textrm{ if } \vert \underline{i} \vert = d-1
\end{array} \right. . $$
For $\underline{i}\in{\Bbb N}^{m}$ and $j\in \{1,\ldots,m\}$, we define the element
$\underline{i}(j)$ of ${\Bbb N}^{m}$ by:
$$i(j):=(i_1,\ldots,i_{j-1},i_j+1,i_{j+1},\ldots,i_m). $$
It suffices to prove:

\begin{claim}\label{clss}
For $\vert i \vert \leq  d-1$, ${\psi}_{\underline{i}}$ is an element of
${\goth g}^{e}$ such that $[\psi _{\underline{i}},u_{j}]+[\psi _{\underline{i}(j)},e]=0$
for $j=1,\ldots,m$.
\end{claim}
Indeed, by Lemma \ref{l4ss2}, if
$$ [\psi _{\underline{i}},u_{j}]+[\psi _{\underline{i}(j)},e]=0 \mbox{ and }
\psi _{\underline{i}} \in \g ^{e} $$
for all $j=1,\ldots,m$, then $\psi _{\underline{i}}\in \z(\g^{e})$.
Furthermore, if
$$ [\psi _{\underline{i}},u_{j}]+[\psi _{\underline{i}(j)},e]=0
\mbox{ and } \psi _{\underline{i}} \in \g ^{e}
\mbox{ and } \psi _{\underline{i}(j)} \in \g ^{e}$$
for all $j=1,\ldots,m$, then $\psi _{\underline{i}}=0$ since
$\z(\g^{e})\cap \g^{f}=0$. So only remains to prove Claim~\ref{clss}.

We prove the claim by induction on $\vert \underline{i} \vert$. Arguing as in the case
$d=1$, we prove the claim for $\vert \underline{i} \vert=0$. We suppose the claim true
for all $\vert \underline{i} \vert\leq l-1$ for some $0 < l \leq d-2$. We have to prove
the statement for all $\vert i \vert\leq l$. By what foregoes and by induction hypothesis,
$\psi _{\underline{i}}= 0$ for $\vert \underline{i} \vert\leq l-2$. For $k=l+1,l+2$, we
consider the ring $\k[\tau _{k}]$ where $\tau _{k}^{k}=0$. Since
$(u_{1}^{*})^{d}$ vanishes on the set of $\k[\tau _{l+1}]$-points
$x=x_0+ x_1 \tau_{l+1} + \cdots+ x_{l}\tau_{l+1}^{l}$ of $Y$ whose source $x_0$ is a
zero of $u_{1}^{*}$,
$$ 0 = [\til{\psi}(e+\tau _{l+1}v),e+\tau _{l+1}v] = \sum_{\vert i \vert = l}
\tau _{l+1}^{l}[\psi _{\underline{i}},e](u_{1}^{*})^{i_{1}}\cdots
(u_{m}^{*})^{i_{m}}(v), $$
for all $v\in{\goth g}^{f}$. So $\psi _{\underline{i}}\in {\goth g}^{e}$ for
$\vert \underline{i}  \vert = l$.

For $\vert \underline{i} \vert$ equal to $l$, the term in
$$\tau _{l+2}^{l+1}(u_{1}^{*})^{i_{1}}\cdots (u_{i_{j-1}}^{*})^{i_{j-1}}
(u_{i_j+1}^{*})^{i_{j}+1} (u_{i_{j+1}}^{*})^{i_{j+1}} \cdots (u_{m}^{*})^{i_{m}}(v) $$
of $[\til{\psi}(e+\tau _{l+2} v ),e+\tau _{l+2}v ]$ is equal to
$[\psi _{\underline{i}(j)},e] + [\psi _{\underline{i}},u_{j}]$. Since
$(u_{1}^{*})^{d}$ vanishes on the set of $\k[\tau _{l+2}]$-points of $Y$ whose source is
a zero of $u_{1}^*$, this term is equal to $0$, whence the claim.
\end{proof}

Recall that $Y'$ is an affine open subset of $Y$ contained in $\Omega $ and having a
nonempty intersection with $p ^{-1}(e)$.

\begin{coro}\label{c2ss2}
Suppose that $e$ has Property $(\mathrm{P})$. Let $\varphi $ be in
$\tk {\k}{\k[Y']}{\goth g}$ such that $\varphi(x)\in[{\goth g},p(x)]$ for all $x \in Y'$.
Then there exists $\psi $ in $\tk {\k}{\k[Y']}{\goth g}$ such that
$\varphi(x)=[\psi(x),p (x)]$ for all $x\in Y'$.
\end{coro}

\begin{proof}
For $a\in\k[Y]$, denote by $D(a)$ the principal open subset defined by $a$.
Let $D(a_1),\ldots,D(a_m)$ be an open covering
of $Y'$ by principal open subsets of $Y$, with $a_1,\ldots,a_k$ in $\k[Y]$. Since
$\varphi$ is a regular map from $Y'$ to ${\goth g}$, there is $m_{i} \geq 0$ such that
$a_{i}^{m_{i}}\varphi$ is the restriction to $Y'$ of some regular map $\varphi_{i}$
from $Y$ to ${\goth g}$. For $m_{i}$  big enough, $\varphi_{i}$ vanishes
on $Y\setminus D(a_i)$; hence $\varphi_{i}(x)\in [{\goth g},p(x)]$ for all $x\in Y$. So,
by Proposition \ref{pss2}, there is a regular map $\psi_{i}$ from $Y'$ to ${\goth g}$ such
that $\varphi_{i}(x)=[\psi_{i}(x),p (x)]$ for all $x\in Y'$. Then for all $x\in D(a_i)$,
we have $\varphi (x)=[a_{i}(x)^{-m_{i}}\psi_{i}(x),p (x)]$. Since $Y'$ is an affine open
subset of $Y$, there exists a regular map $\psi$ from $Y'$ to ${\goth g}$ which satisfies
the condition of the corollary.
\end{proof}

We are now in position to prove the main result of this section:

\begin{theo}\label{tss2}
The equality $\ind\g^e=\rg$ holds if and only if $e$ has Property $(\mathrm{P})$.
\end{theo}

\begin{proof}
By Corollary \ref{c2ss2}, if $e$ has Property $(\mathrm{P})$, then Assertion (B)
of Theorem \ref{tss1} is satisfied. Conversely, suppose that $\ind\g^e=\rg$ and show that
$e$ has Property (P). By Theorem~\ref{tss1}, (A)$\Rightarrow$(B), Assertion (B)
is satisfied. We choose an affine open subset $Y'$ of $Y$, contained in $\Omega $,
such that $Y'\cap p^{-1}(e)\not= \empty$ and verifying the condition of the assertion (B).
Let $z\in\z(\g^e)$.
Recall that the map $\varphi_{z}$ is  defined by~(\ref{eq:phi-z}).
Let $x$ be in $Y'$. If $u_{1}^{*}(x)\not=0$, then
$\varphi_{z}(x)$ belongs to $[{\goth g},p(x)]$ by~(\ref{eq:phi-z}). If
$u_{1}^{*}(x)=0$ , then by Lemma~\ref{l4ss2}, $\varphi_{z}(x)$ belongs to
$[e,{\goth g}]$. So there exists a regular map $\psi $ from $Y'$ to ${\goth g}$
such that $\varphi_{z}(x)=[\psi(x),p(x)]$ for all $x\in Y'$ by Assertion (B).
Hence we have
$$[z-u_{1}^{*}\psi (x),p(x)] = 0 ,$$
for all $x\in Y'$ since $(u_{1}^{*}\varphi_{z})(x)=[z,p(x)]$ for all $x\in Y$. So
$\alpha (x)$ contains $z$ for all $x$ in $\Omega  \cap Y' \cap p ^{-1}(e)$.
Since $p^{-1}(e)$ is irreducible, we deduce that $e$ has Property (P).
\end{proof}

\subsection{A new formulation of the property (P)}\label{ss3}
Recall that Property (P) is introduced in Definition~\ref{dss1}. As has been
noticed in the proof of Lemma~\ref{lss2}, the morphism
$G\times\slo \to \g, (g,x)\mapsto g(x)$ is smooth. As a consequence, the set $\slo _{\r}$
of $v \in \slo$ such that $v$ is regular is a nonempty open subset of $\slo$.
For $x$ in $\slo_{\r}$, ${\goth g}^{e+t(x-e)}$ has dimension $\rg$ for all $t$ in a
nonempty open subset of $\k$ since $x=e+(x-e)$ is regular. Furthermore, since $\k$ has
dimension $1$, \cite[Ch.~VI, Theorem 1]{Sh} asserts that there is a unique regular map
$$ \beta _{x}: \k\to \mathrm{Gr}_{\rg}(\g)$$
satisfying $\beta _{x}(t)={\goth g}^{e+t(x-e)}$ for  all $t$ in a nonempty open subset of
$\k$.

Recall that $Y$ is an affine open subset of $\blw$ with $Y\cap p^{-1}(e)\not=\emptyset$
and that $u_{1}^{*},v_{2}^{*},\ldots,v_{m}^{*}$ is a global coordinates system of $Y$,
cf.~(\ref{eq:vi}). Let $\slo'_{\r}$ be the subset of $x$ in $\slo_{\r}$ such
that $u_{1}^{*}(x)\not=0$. For $x$ in $\slo_{\r}'$, we denote by $\til{x}$ the element
of $Y$ whose coordinates are $0, v_{2}^{*}(x),\ldots, v_{m}^{*}(x)$.

\begin{lemma}\label{lss3}
Let $x$ be in $\slo'_{\r}$.

{\rm (i)} The subspace $\beta _{x}(0)$ is contained in ${\goth g}^{e}$.

{\rm (ii)} If $\til{x}\in\Omega $, then $\alpha (\til{x})=\beta _{x}(0)$.
\end{lemma}

\begin{proof}
(i) The map $\beta _{x}$ is a regular map and $[\beta _{x}(t),e+t(x-e)]=0$ for all $t$
in a nonempty open subset of $\k$. So, $\beta _{x}(0)$ is contained in ${\goth g}^{e}$.

(ii) Since $\slo'_{\r}$ has an empty intersection with the nullvariety of $u_{1}^{*}$
in $\slo$, the restriction of $p$ to $p^{-1}(\slo'_{\r})$ is an isomorphism from
$p^{-1}(\slo'_{\r})$ to $\slo'_{\r}$. Furthermore,
$\beta _{x}(t)=\alpha (p^{-1}(e+tx-te))$ for any $t$ in $\k$ such that $e+t(x-e)$
belongs to $\slo'_{\r}$ and $p^{-1}(e+tx-te)$ goes to $\til{x}$ when $t$ goes to $0$.
Hence $\beta _{x}(0)$ is equal to $\alpha (\til{x})$ since $\alpha $ and $\beta $ are
regular maps.
\end{proof}

\begin{coro}\label{css3}
The element $e$ has Property {\rm (P)} if and only if $\z(\g^e) \subset \beta _{x}(0)$
for all $x$ in a nonempty open subset of $\slo _{\r}$.
\end{coro}

\begin{proof}
The map $x\mapsto \til{x}$ from ${\slo}'_{\r}$ to $Y$ is well-defined and its image is an
open subset of $Y\cap p ^{-1}(e)$. Let ${\slo}''_{\r}$ be the set of $x\in {\slo}'_{\r}$
such that $\til{x}\in \Omega $ and let $Y''$ be the image of ${\slo}''_{\r}$ by the map
$x\mapsto \til{x}$. Then ${\slo}''_{\r}$ is open in $\slo_{\r}$ and $Y''$ is dense in
$\Omega \cap p ^{-1}(e)$ since $p^{-1}(e)$ is irreducible. Furthermore, the image of
a dense open subset of $\slo ''_{\r}$ by the map $x\mapsto \til{x}$ is dense in $Y''$.
Since $\alpha $ is regular, $e$ has property (P) if and only if $\alpha (x)$ contains
$\z(\g^e)$ for all $x$ in a dense subset of $Y''$. By Lemma~\ref{lss3}(ii), the latter
property is equivalent to the fact that $\beta _{x}(0)$ contains $\z(\g^e)$ for all $x$
in a dense open subset of $\slo''_{\r}$.
\end{proof}

\begin{coro}\label{c2ss3}
{\rm (i)} If $\z(\g^e)$ is generated by $\varphi_1(e),\ldots,\varphi_{\rg}(e)$, then $e$
has Property {\rm (P)}.

{\rm (ii)} If $\z(\g^e)$ has dimension 1, then $e$ has Property {\rm (P)}.
\end{coro}

\begin{proof}
Recall that $\varphi _{i}(e)$ belongs to $\z(\g^{e})$, for all $i=1,\ldots,\rg$,  by
Lemma~\ref{ldi}(i). Moreover,  for all $x$ in $\slo_{\r}$ and all
$i=1,\ldots,\rg$, $\varphi _{i}(e+t(x-e))$ belongs to ${\goth g}^{e+t(x-e)}$ for any $t$
in $\k$. So by continuity, $\varphi _{i}(e)$ belongs to $\beta _{x}(0)$.
As a consequence, whenever $\z(\g^{e})$ is generated by
$\poi e{}{,\ldots,}{}{\varphi }{1}{\rg}$, $e$ has Property {\rm (P)}
by Corollary~\ref{css3}.

(ii) is an immediate consequence of (i) since $\varphi_1(e)=e$ by our choice of $d_1$.
\end{proof}

\section{Proof of Theorem~\ref{tint2} for rigid nilpotent orbits}\label{ar}
We intend to prove in this section the following theorem:

\begin{theo}\label{tar}
Suppose that $\g$ is reductive and let $e$ be a rigid nilpotent element of $\g$.
Then the index of ${\goth g}^{e}$ is equal to $\rg$.
\end{theo}

Theorem~\ref{tar} will complete the proof of Theorem~\ref{tint2} by Theorem~\ref{tcr}.
As explained in introduction, we can assume that $\g$ is simple. We consider two cases,
according to $\g$ has classical type or exceptional type.

\subsection{The classical case}\label{ar1}
Assume that ${\goth g}$ is simple of classical type. More precisely, assume that
${\goth g}$ is one of the Lie algebras ${\goth {sl}}_{\rg + 1}(\k)$,
${\goth {so}}_{2\rg + 1}(\k)$, ${\goth {sp}}_{2\rg }(\k)$, ${\goth {so}}_{2\rg }(\k)$.

\begin{lemma}\label{lrss3}
Let $m$ be a positive integer such that $x^{m}-\tr x^{m}$ belongs to ${\goth g}$ for
all $x$ in ${\goth g}$. Then $e^{m}$ belongs to the subspace generated by
$\poi e{}{,\ldots,}{}{\varphi }{1}{\rg}$.
\end{lemma}

\begin{proof}
Recall that $L_{{\goth g}}$ is the submodule of elements $\varphi $ of
$\tk {\k}{\e Sg}{\goth g}$ such that $[x,\varphi (x)]=0$ for all $x$ in ${\goth g}$.
According to \cite{Di3}, $L_{{\goth g}}$ is a free module generated by the
$\varphi _{i}'s$. For all $x$ in ${\goth g}$, $[x,x^{m}]=0$. Hence there exist
polynomial functions $\poi a1{,\ldots,}{\rg}{}{}{}$ on ${\goth g}$ such that
$$ x^{m}-\tr x^{m} = a_{1}(x)\varphi _{1}(x) + \cdots + a_{\rg}(x)\varphi _{\rg}(x)$$
for all $x$ in ${\goth g}$, whence the lemma.
\end{proof}

\begin{theo}\label{tar1}
Let $e$ be a rigid nilpotent element. Then $\z(\g^{e})$ is generated by powers of $e$. In
particular, the index of $\g^{e}$ is equal to $\rg$.
\end{theo}

\begin{proof}
Let us prove the first assertion. If $\g$ has type A or C, then $\z(\g^{e})$ is generated
by powers of $e$ by~\cite[Th\'eor\`eme 1.1.8]{Mo3} or~\cite{Y2}. So we can assume that
$\g$ has type B or D.

Set $n:=2\rg+1$ if ${\goth g}$ has type ${\mathrm {B}}_{\rg}$ and $n:=2\rg$ if
${\goth g}$ has type ${\mathrm {D}}_{\rg}$. Denote by $(\poi n1{,\ldots,}{k}{}{}{})$,
with $\poi n1{ \geq \cdots \geq }{k}{}{}{}$, the partition of $n$ corresponding to the
nilpotent element $e$. By~\cite[Th\'eor\`eme 1.1.8]{Mo3} or~\cite{Y2},
$\z({\goth g}^{e})$ is \emph{not} generated by powers of $e$ if and only if $n_{1}$ and
$n_{2}$ are both odd integers and $n_{3}<n_{2}$. On the other hand, since $e$ is rigid,
$n_{k}$ is equal to $1$, $n_{i}\leq n_{i+1}\leq n_{i}+1$ and all odd integers of the
partition $(\poi n1{,\ldots,}{k}{}{}{})$ have a multiplicity different from $2$
\cite[ch.~II]{Ke,Sp} or~\cite[Corollary 7.3.5]{CMa}. Hence, the preceding criterion is
not satisfied for $e$. Then, the second assertion results from Lemma \ref{lrss3},
Corollary~\ref{c2ss3}(i) and Theorem~\ref{tss2}.
\end{proof}

\begin{rema}\label{rar1}
Yakimova's proof of Elashvili's conjecture in the classical case is shorter and more
elementary~\cite{Y1}. The results of Section~\ref{ss} will serve the exceptional case
in a more relevant way.
\end{rema}

\subsection{The exceptional case}\label{ar2}
We let in this subsection $\g$ be simple of exceptional type and we assume
that $e$ is a nonzero rigid nilpotent element of $\g$.
The dimension of the center of centralizers
of nilpotent elements has been recently described in \cite[Theorem 4]{LaTe}.
On the other hand, we have explicit computations for the rigid nilpotent orbits
in the exceptional types due to A.G.~Elashvili.
These computations are collected
in~\cite[Appendix of Chap.~II]{Sp} and a complete version was published later
in~\cite{E2}.
From all this, we observe that the center of $\g^e$ has dimension $1$ in
most cases.
In more details, we have:

\begin{prop}\label{par2}
Let $e$ be nonzero rigid nilpotent element of $\g$.

{\rm (i)} Suppose that $\g$ has type $\mathrm{G}_2$, $\mathrm{F}_4$ or $\mathrm{E}_6$.
Then $\dim\z(\g^e)=1$.

{\rm (ii)} Suppose that $\g$ has type $\mathrm{E}_7$. If $\g^e$ has dimension $41$, then
$\dim \z(\g^e)=2$; otherwise $\dim \z(\g^e)=1$.

{\rm (iii)} Suppose that $\g$ has type $\mathrm{E}_8$. If $\g^e$ has dimension $112$,
$84$, $76$, or $46$, then $\dim \z(\g^e)=2$, if $\g^e$ has dimension $72$, then
$\dim \z(\g^e)=3$; otherwise $\dim \z(\g^e)=1$.
\end{prop}

By Corollary~\ref{c2ss3}(ii), $\ind\g^e=\rg$ whenever $\dim \z(\g^e)=1$.
So, as an immediate consequence of Proposition~\ref{par2}, we obtain:

\begin{coro}\label{car2}
Suppose that either $\g$ has type $\mathrm{G}_{2}$, $\mathrm{F}_{4}$, $\mathrm{E}_{6}$,
or $\g$ has type $\mathrm{E}_{7}$ and $\dim \g^{e} \neq 41$, or $\g$ has type
$\mathrm{E}_{8}$ and \hbox{$\dim \g^{e}\not\in \{112,84,76,72,46\}$}.
Then $\dim \z(\g^e)=1$ and the index of $\g ^{e}$ is equal to $\rg$.
\end{coro}

According to Corollary~\ref{car2}, it remains 7 cases;
there are indeed two rigid nilpotent orbits of codimension 46 in E$_8$.
We handle now these remaining cases.
We process here in a different way; we study technical conditions on $\g^e$
under which $\ind\g^e=\rg$. For the moment, we state general results about the index.

Let ${\goth a}$ be an algebraic Lie algebra. Recall that the stabilizer of
$\xi \in \a^*$ for the coadjoint representation is denoted by ${\goth a}^{\xi}$ and
that $\xi $ is regular if $\dim {\goth a}^{\xi }=\ind \a$.
Choose a commutative subalgebra ${\goth t}$ of ${\goth a}$ consisted of semisimple
elements of $\a$ and denote by $\z_{\a}({\goth t})$ the centralizer of ${\goth t}$ in
${\goth a}$. Then
${\goth a}={\z_{\a}({\goth t})} \oplus [{\goth t},{\goth a}]$.
The dual $\z_{\a}({\goth t})^{*}$ of $\z_{\a}({\goth t})$ identifies to the orthogonal
complement of $[{\goth t},{\goth a}]$ in ${\goth a}^{*}$. Thus,
$\xi \in \z_{\a}({\goth t})^*$ if and only if ${\goth t}$ is contained in $\a^{\xi}$.

\begin{lemma}\label{lar2}
Suppose that there exists $\xi $ in ${\goth z}_{{\goth a}}({\goth t})^{*}$ such that
$\dim ({\goth a}^{\xi }\cap [{\goth t},{\goth a}]) \leq 2$.
Then
$$\ind {\goth a} \leq \ind {\goth z}_{{\goth a}}({\goth t})+ 1.$$
\end{lemma}

\begin{proof}
Let $T$ be the closure in
${\goth z}_{{\goth a}}({\goth t})^{*}\times \ec {Gr}{}{[{\goth t},{\goth a}]}{}3$ of the
subset of elements $(\eta ,E)$ such that $\eta $ is a regular element of
${\goth z}_{{\goth a}}({\goth t})^{*}$ and $E$ is contained in ${\goth a}^{\eta}$.
The image $T_1$ of $T$ by the projection from
${\goth z}_{{\goth a}}({\goth t})^{*}\times \ec {Gr}{}{[{\goth t},{\goth a}]}{}3$
to ${\goth z}_{{\goth a}}({\goth t})^{*}$ is
closed in ${\goth z}_{{\goth a}}({\goth t})^{*}$.
By hypothesis, $T_1$ is not equal to ${\goth z}_{{\goth a}}({\goth t})^{*}$ since
for all $\eta$ in $T_1$, $\dim ({\goth a}^{\eta }\cap [{\goth t},{\goth a}]) \geq 3$.
Hence there exists a regular element $\xi_0 $ in ${\goth z}_{{\goth a}}({\goth t})^{*}$
such that $\dim ({\goth a}^{\xi_0 }\cap [{\goth t},{\goth a}]) \leq 2$.
Since ${\goth t}$ is contained in ${\goth a}^{\xi_0 }$,
$$ {\goth a}^{\xi_0 } = {\goth z}_{{\goth a}}({\goth t})^{\xi_0 }\oplus
{\goth a}^{\xi_0 }\cap [{\goth t},{\goth a}] .$$
If $[{\goth t},{\goth a}]\cap {\goth a}^{\xi_0 }=\{0\}$ then $\ind {\goth a}$ is
at most $\ind {\goth z}_{{\goth a}}({\goth t})$. Otherwise, ${\goth a}^{\xi_0 }$ is not a
commutative subalgebra since ${\goth t}$ is contained in ${\goth a}^{\xi_0 }$.
Hence $\xi_0 $ is not a regular element of ${\goth a}^{*}$, so
$\ind {\goth a} < \dim {\goth a}^{\xi_0 }$. Since
$\dim {\goth a}^{\xi_0 }\leq \ind {\goth z}_{{\goth a}}({\goth t})+2$, the lemma follows.
\end{proof}

From now on, we assume that ${\goth a}=\g ^{e}$. As a rigid nilpotent element of $\g$,
$e$ is a nondistinguished nilpotent element. So we can choose a nonzero commutative subalgebra ${\goth t}$
of $\g^{e}$ consisted of semisimple elements. Denote by $\l$ the centralizer of
${\goth t}$ in ${\goth g}$. As a Levi subalgebra of $\g$, $\l$ is a reductive Lie algebra
whose rank is $\rg$. Moreover its dimension is strictly smaller than $\dim \g$.
In  the preceding notations, we have
$ {\goth z}_{\g^e}({\goth t}) = {\goth z}_{\g}({\goth t})^e=\l^e$.
Let ${\goth t}_{1}$ be a commutative subalgebra of ${\goth l}^{e}$ containing ${\goth t}$
and consisting of semisimple elements of $\l$. Then $[{\goth t},{\goth g}^{e}]$ is stable
under the adjoint action of ${\goth t}_{1}$. For $\lambda $ in ${\goth t}_{1}^{*}$,
denote by $\g ^{e}_{\lambda }$ the $\lambda$-weight space of the adjoint action of
${\goth t}_{1}$ in ${\goth g}^{e}$.

\begin{lemma}\label{l2ar2}
Let $\lambda \in {\goth t}_{1}^{*}$ be a nonzero weight of the adjoint action of
${\goth t}_{1}$ in $\g^{e}$. Then $-\lambda $ is also a weight for this action and
$\lambda$ and $-\lambda$ have the same multiplicity. Moreover, $\g ^{e}_{\lambda }$ is
contained in $[{\goth t},\g^{e}]$ if and only if the restriction of $\lambda $ to
${\goth t}$ is not identically zero.
\end{lemma}

\begin{proof}
By definition, $\g ^{e}_{\lambda }\cap {\goth l}^{e}=\{0\}$ if and only if
the restriction of $\lambda $ to ${\goth t}$ is not identically zero.
So ${\goth g}^{e}_{\lambda }$ is contained in $[{\goth t},{\goth g}^{e}]$ if and only if
the restriction of $\lambda $ to ${\goth t}$ is not equal to $0$ since
$$ {\goth g}^{e}_{\lambda } = ({\goth g}^{e}_{\lambda }\cap {\goth l}^{e}) \oplus
({\goth g}^{e}_{\lambda }\cap [{\goth t},{\goth g}^{e}]) .$$
The subalgebra ${\goth t}_{1}$ is contained in a reductive factor of ${\goth g}^{e}$. So
we can choose $h$ and $f$ such that ${\goth t}_{1}$ is contained in
${\goth g}^{e}\cap {\goth g}^{f}$. As a consequence, any weight of the adjoint action of
${\goth t}_{1}$ in ${\goth g}^{f}$ is a weight of the adjoint action of ${\goth t}_{1}$
in ${\goth g}^{e}$ with the same multiplicity. Furthermore, the ${\goth t}_{1}$-module
${\goth g}^{f}$ for the ajoint action is isomorphic to the ${\goth t}_{1}$-module
$({\goth g}^{e})^{*}$ for the coadjoint action. So $-\lambda $ is a weight of the adjoint
action of ${\goth t}_{1}$ in ${\goth g}^{f}$ with the same multiplicity as $\lambda $.
Hence $-\lambda $ is a weight of the adjoint action of ${\goth t}_{1}$ in
${\goth g}^{e}$ with the same multiplicity as $\lambda $, whence the lemma.
\end{proof}

Choose pairwise different elements $\poi {\lambda }1{,\ldots,}{r}{}{}{}$ of
${\goth t}_{1}^{*}$ so that the weights of the adjoint action of ${\goth t}_{1}$ in
$\g^{e}$ which are not identically zero on ${\goth t}$ are precisely the elements
$\pm \lambda _{i}$. For $i=1,\ldots,r$, let
$\poi v{i,1}{,\ldots,}{i,m_{i}}{}{}{}$ and $\poi w{i,1}{,\ldots,}{i,m_{i}}{}{}{}$ be
basis of $\g ^{e}_{\lambda _{i}}$ and $\g ^{e}_{-\lambda _{i}}$ respectively.
Then we set:
$$ q_{i} := \det \left (
[v_{i,k},w_{i,l}] \right )_{1 \leq k,l \leq m_i} \in \mathrm{S}(\l^e).$$

\begin{prop}\label{p2ar2}
Suppose that $\ind {\l^e}=\rg$ and suppose that one of the following two conditions is
satisfied:
\begin{enumerate}
\item   for $i=1,\ldots,r$, $q_{i}\neq 0$,
\item   there exists $j$ in $\{1,\ldots,r\}$ such that $q_{i}\neq 0$ for all
$i\neq j$ and such that the basis $\poi v{j,1}{,\ldots,}{j,m_{j}}{}{}{}$ and
$\poi w{j,1}{,\ldots,}{j,m_{j}}{}{}{}$ can be chosen so that
$$\det \left ([v_{j,k},w_{j,l}] \right )_{1 \leq k,l \leq m_j-1} \not=0  .$$
\end{enumerate}
Then, $\ind \g ^{e}=\rg$.
\end{prop}

\begin{proof}
First, observe that $\ind \g^e -\ind \g$ is an even integer.
Indeed, we have:
\begin{eqnarray*}
&& \ind \g^e - \ind \g= (\ind \g^e - \dim \g^e) + (\dim \g^e - \dim\g) +
(\dim \g -\ind \g).
\end{eqnarray*}
But the integers $\ind \g^e - \dim \g^e$, $\dim \g^e - \dim\g$ and $\dim \g -\ind \g$ are
all even integers. Thereby, if $\ind \g^e \leq \ind\g +1$, then $\ind \g^e \leq \ind \g$.
In turn, by Vinberg's inequality (cf.~Introduction), we have $\ind \g^e \geq \ind \g$.
Hence, it suffices to prove $\ind \g^e \leq \ind \l^e +1$ since our hypothesis says that
$\ind \l^e=\rg=\ind\g$. Now, by Lemma \ref{lar2}, if there exists  $\xi $ in
$(\l ^{e})^{*}$ such that
$(\g ^{e})^{\xi }\cap [{\goth t},\g ^{e}]$ has dimension at most $2$, then we are done.

Denote by $\l _{1}$ the centralizer of ${\goth t}_{1}$ in $\g $. Then $\l _{1}$ is
contained in $\l $ and $\l ^{e}= \l _{1}^{e}\oplus [{\goth t}_{1},\l^{e}]$ and
 $(\l _{1}^{e})^{*}$ identifies to the orthogonal of $[{\goth t}_{1},\l^{e}]$ in the dual
of $\l ^{e}$. Moreover, for $i=1,\ldots,r$, $q_{i}$ belongs to $\es S{\l _{1}^{e}}$.
For $\xi $ in $({\goth l}_{1}^{e})^*$, denote by $B_{\xi}$ the bilinear form
$$\begin{array}{ccc}
 [{\goth t},{\goth g}^{e}]\times [{\goth t},{\goth g}^{e}]& \longrightarrow& \k \\
(v,w) &\longmapsto &\xi([v,w])
\end{array}$$
and denote by $\ker B_{\xi}$ its kernel.
For $i=1,\ldots,r$, $-q_{i}(\xi )^{2}$ is the determinant of the restriction of $B_{\xi }$
to the subspace
$$ ({\goth g}^{e}_{\lambda _{i}}\oplus {\goth g}^{e}_{-\lambda _{i}}) \times
({\goth g}^{e}_{\lambda _{i}}\oplus {\goth g}^{e}_{-\lambda _{i}}) $$
in the basis $\poi v{i,1}{,\ldots,}{i,m_{i}}{}{}{},\poi w{i,1}{,\ldots,}{i,m_{i}}{}{}{}$.

If (1) holds, we can find $\xi $ in $(\l _{1}^{e})^{*}$ such that $\ker B_{\xi }$ is
zero. If (2) holds, we can find $\xi $ in $(\l _{1}^{e})^{*}$ such that $\ker B_{\xi }$
has dimension $2$ since $B_{\xi }$ is invariant under the adjoint action of
${\goth t}_{1}$. But $\ker B_{\xi }$ is equal to
$(\g ^{e})^{\xi }\cap [{\goth t},\g ^{e}]$. Hence such a $\xi $ satisfies the
required inequality and the proposition follows.
\end{proof}

The proof of the following proposition is given in
Appendix~\ref{bp} since it relies on explicit computations:

\begin{prop}\label{p3ar2}
{\rm (i)} Suppose that either $\g$ has type ${\mathrm {E}}_{7}$ and $\dim \g^{e}=41$ or,
$\g$ has type ${\mathrm {E}}_{8}$ and $\dim \g^{e}\in \{112,72\}$.
Then, for suitable choices of ${\goth t}$ and ${\goth t}_{1}$, Condition {\rm (1)} of
Proposition~{\rm \ref{p2ar2}} is satisfied.

{\rm (ii)} Suppose that $\g$ has type ${\mathrm {E}}_{8}$ and  that $\g^{e}$
has dimension $84$, $76$, or $46$. Then, for  suitable choices of ${\goth t}$ and
${\goth t}_{1}$, Condition {\rm (2)} of Proposition~{\rm \ref{par2}} is satisfied.
\end{prop}

\subsection{Proof of Theorem~\ref{tint2}}
We are now in position to complete the proof of Theorem~\ref{tint2}:

\begin{proof}
[Proof of Theorem~\ref{tint2}]
We argue by induction on the dimension of $\g$.
If $\g$ has dimension 3, the statement is known.
Assume now that
$\ind \l^{e'}=\rk \l$ for any reductive Lie algebras $\l$ of dimension at most
$\dim \g - 1$ and any $e' \in \nil{\l}$.
Let $e\in \nil{\g}$ be a nilpotent element of $\g$.
By Theorem~\ref{tcr} and Theorem~\ref{tar1}, we can assume that $e$ is rigid and that
$\g$ is simple of exceptional type. Furthermore by Corollary~\ref{car2}, we can assume
that $\dim\z(\g^e)>1$. Then we consider the different cases given by
Proposition~\ref{p3ar2}.

If, either $\g$ has type E$_{7}$ and $\dim \g^{e}= 41$, or $\g$ has type
E$_{8}$ and $\dim \g^{e}$ equals 112, 72, or 46, then Condition (1) of
Proposition~\ref{p2ar2} applies for suitable choices of
${\goth t}$ and ${\goth t}_{1}$ by Proposition~\ref{p3ar2}.
Moreover, if $\l=\z_{\g}({\goth t})$, then
$\l$ is a reductive Lie algebra of rank $\rg$ and strictly contained in ${\goth g}$.
So, from our induction hypothesis, we deduce that $\ind \g^e=\rg$ by
Proposition~\ref{p2ar2}.

If $\g$ has type E$_{8}$ and $\dim \g^{e}$ equals 84, 76, or 46, then
Condition (2) of Proposition~\ref{p2ar2} applies for suitable choices of
${\goth t}$ and ${\goth t}_{1}$ by Proposition~\ref{p3ar2}.
Arguing as above, we deduce
that $\ind \g^e=\rg$.
\end{proof}

\appendix \label{A}

\section{Proof of Proposition~\ref{p3ar2}: explicit computations.}\label{bp}
This appendix aims to prove Proposition~\ref{p3ar2}.
We prove Proposition~\ref{p3ar2} for each case
by using explicit computations made with the help of \texttt{GAP};
our programmes are presented below
(two cases are detailed; the other ones are similar).
Explain the general approach. In our programmes, $x[1],\ldots$ are root vectors generating
the nilradical of the Borel subalgebra ${\goth b}$ of ${\goth g}$ and the representative
$e$ (denoted by \texttt{e} in our programmes) of the rigid orbit is chosen so that $e$ and
$h$ belong to ${\goth b}$ and ${\goth h}$ respectively. The element $e$ is given by the
tables of~\cite{GQT}. In fact, in~\cite{GQT}, they use the programme Lie which induces
minor changes in the numbering. Then, we exhibit suitable tori ${\goth t}$ and
${\goth t_1}$ of $\g$ contained in $\g^e$ which satisfies conditions (1)
or (2) of Proposition~\ref{p2ar2}. In each case, our torus ${\goth t}$ is one
dimensional;  we define it by a generator, called \texttt{t} in our programmes.
Its centralizer in $\g^e$ is denoted by \texttt{le}.
The torus ${\goth t}_{1}$ has dimension at most $4$.
It is defined by a basis denoted by \texttt{Bt1}.
The weights of ${\goth t}_{1}$ for the adjoint action of ${\goth t}_{1}$
on $\g^e$ are given by their values on the basis \texttt{Bt1} of ${\goth t}_{1}$.
We list in a matrix \texttt{W} almost all weights which have a positive value at
\texttt{Bt1}.
The other weights have multiplicity $1$.
In our programmes, by the term
\texttt{S} we check that no weight is forgotten; this term has to be zero.
Then, the matrices corresponding to the weights given by \texttt{W} are given by a
function \texttt{A}.
Their determinants correspond to the $q_{i}$'s in the notations
of Proposition~\ref{p2ar2}.
If there is only one other weight, the corresponding matrix is denoted by \texttt{a}.
At last, we verify that these matrices have the desired property depending
on the situations (i) or (ii) of Proposition~\ref{p3ar2}.

As examples, we detail below two cases:

(1) the case of E$_7$ with $\dim \g^e=41$
where we intend to check that Condition (1) of Proposition~\ref{p2ar2} is satisfied;

(2) the case of E$_8$ with $\dim \g^e=84$
where we intend to check that Condition (2) of Proposition~\ref{p2ar2} is satisfied.\\

(1) E$_{7}$, $\dim {\goth g}^{e}=41$:
In this case, with our choices, $\dim {\goth t}=1$, $\dim \l^e=23$ and
$\dim {\goth t}_{1}=3$. The order of matrices to be considered is at most $2$.

\begin{verbatim}
L := SimpleLieAlgebra("E",7,Rationals);;  R := RootSystem(L);;
x := PositiveRootVectors(R);; y := NegativeRootVectors(R);;
e := x[14]+x[26]+x[28]+x[49];;
c := LieCentralizer(L,Subspace(L,[e]));Bc := BasisVectors(Basis(c));;
> <Lie algebra of dimension 41 over Rationals>
z := LieCentre(c);; Bz := BasisVectors(Basis(z));;
t := Bc[Dimension(c)];;
le := LieCentralizer(L,Subspace(L,[t,e]));
> <Lie algebra of dimension 23 over Rationals>
n := function(k)
     if k=2 then return 1;;
     elif k=-2 then return 1;;
     elif k=1 then return 8;;
     elif k=-1 then return 8;; fi;; end;;
#The function n assigns to each weight of t the dimension of the corresponding
#weight subspace.
M := function(k)  local m;;
     m := function(j,k)
          if j=1 then return Position(List([1..Dimension(c)],
          i->t*Bc[i]-k*Bc[i]),0*x[1]);;
          else return m(j-1,k)+Position(List([m(j-1,k)+1..Dimension(c)],
          i->t*Bc[i]-k*Bc[i]),0*x[1]);;
          fi;;
     end;;
     return List([1..n(k)],i->m(i,k));;
end;;
Bt1 := [Bc[41],Bc[40],Bc[39]];;
N := function(k,p) local n;;
     n := function(j,k,p)
          if j=1 then return Position(List([1..8],
          i->Bt1[2]*Bc[M(k)[i]]-p*Bc[M(k)[i]]),0*x[1]);;
          else return n(j-1,k,p)+Position(List([n(j-1,k,p)+1..8],
          i->Bt1[2]*Bc[M(k)[i]]-p*Bc[M(k)[i]]),0*x[1]);;
          fi;;
     end;;
     return List([1..4],i->M(k)[n(i,k,p)]);;
end;;
r := function(t)
     if t=1 then return 1;
     elif t=-1 then return 1;;
     elif t=0 then return 2;;
     fi;;
end;;
Q := function(k,s,t) local q;;
     q := function(j,k,s,t)
          if j=1 then return Position(List([1..4],
          i->Bt1[3]*Bc[N(k,s)[i]]-t*Bc[N(k,s)[i]]),0*x[1]);;
          else return q(j-1,k,s,t)+Position(List([q(j-1,k,s,t)+1..4],
          i->Bt1[3]*Bc[N(k,s)[i]]-t*Bc[N(k,s)[i]]),0*x[1]);;
          fi;;
     end;;
     return List([1..r(t)],i->N(k,s)[q(i,k,s,t)]);;
end;;
W := [[1,1,1],[1,-1,1],[1,1,-1],[1,-1,-1],[1,1,0],[1,-1,0]];;
S := 2*(1+Sum(List([1..Length(W)],i->Length(Q(W[i][1],W[i][2],W[i][3])))))
     +Dimension(le)-Dimension(c);
> 0
A := function(i) return List([1..r(W[i][3])],t->List([1..r(W[i][3])],
     s->Bc[Q(W[i][1],W[i][2],W[i][3])[s]]*Bc[Q(-W[i][1],-W[i][2],-W[i][3])[t]]));;
end;;
A(1);A(2);A(3);A(4);A(5);A(6);
> [ [ (-1)*v.63 ] ]
> [ [ v.63 ] ]
> [ [ v.63 ] ]
> [ [ (-1)*v.63 ] ]
> [ [ (-1)*v.57+(-1)*v.60, (-1)*v.63 ], [ (-1)*v.63, 0*v.1 ] ]
> [ [ (-1)*v.57+(-1)*v.60, (-1)*v.63 ], [ (-1)*v.63, 0*v.1 ] ]
a := Bc[M(2)[1]]*Bc[M(-2)[1]];
> v.133
\end{verbatim}
In conclusion, Condition (1) of Proposition~\ref{p2ar2} is satisfied
for ${\goth t}:=\k$\texttt{t} and
${\goth t_1}:=$span(\texttt{Bt1}).\\

(2) E$_{8}$, $\dim {\goth g}^{e}=84$:
In this case, with our choices, $\dim {\goth t}=1$, $\dim \l^e=48$ and
$\dim {\goth t}_{1}=3$. The matrix \texttt{A(7)} has order $5$ and it is
singular of rank $4$. The order of the other matrices is at most $2$.

\begin{verbatim}
L := SimpleLieAlgebra("E",8,Rationals);; R := RootSystem(L);;
x := PositiveRootVectors(R);; y := NegativeRootVectors(R);;
e := x[54]+x[61]+x[77]+x[97];;
c := LieCentralizer(L,Subspace(L,[e])); Bc := BasisVectors(Basis(c));;
> <Lie algebra of dimension 84 over Rationals>
z := LieCentre(c);; Bz := BasisVectors(Basis(z));;
t := Bc[Dimension(c)];;
le := LieCentralizer(L,Subspace(L,[t,e]));
> <Lie algebra of dimension 48 over Rationals>
n := function(k)
     if k=2 then return 1;;
     elif k=-2 then return 1;;
     elif k=1 then return 17;;
     elif k=-1 then return 17;;
     fi;;
end;;
M := function(k)  local m;;
     m := function(j,k)
          if j=1 then return Position(List([1..Dimension(c)],
          i->Bc[84]*Bc[i]-k*Bc[i]),0*x[1]);;
          else return m(j-1,k)+Position(List([m(j-1,k)+1..Dimension(c)],
          i->Bc[84]*Bc[i]-k*Bc[i]), 0*x[1]);;
          fi;;
     end;;
     return List([1..n(k)],i->m(i,k));;
end;;
r := function(k,t)
     if k=1 and t=1  then return 4;;
     elif k=-1 and t=-1  then return 4;;
     elif k=1 and t=-1  then return 4;;
     elif k=-1 and t=1  then return 4;;
     elif k=1 and t=0  then return 9;;
     elif k=-1 and t=0  then return 9;;
     fi;;
end;;
Bt1 := [Bc[84],Bc[83],Bc[82]];;
N := function(k,t) local p;;
     p := function(j,k,t)
          if j=1 then return Position(List([1..n(k)],
          i->Bt1[2]*Bc[M(k)[i]]-t*Bc[M(k)[i]]),0*x[1]);;
          else return p(j-1,k,t)+Position(List([p(j-1,k,t)+1..n(k)],
          i->Bt1[2]*Bc[M(k)[i]]-t*Bc[M(k)[i]]),0*x[1]);;
          fi;;
     end;;
     return List([1..r(k,t)],i->M(k)[p(i,k,t)]);;
end;;
m := function(k,s,t)
     if k=1 and s=1 and t=-1 then return 2;;
     elif k=-1 and s=-1 and t=1 then  return 2;;
     elif k=1 and s=1 and t=0 then return 2;;
     elif k=-1 and s=-1 and t=0 then return 2;;
     elif k=1 and s=-1 and t=1 then  return 2;;
     elif k=-1 and s=1 and t=-1 then return 2;;
     elif k=1 and s=-1 and t=0 then return 2;;
     elif k=-1 and s=1 and t=0 then  return 2;;
     elif k=1 and s=0 and t=1 then return 2;;
     elif k=-1 and s=0 and t=-1 then return 2;;
     elif k=1 and s=0 and t=-1 then return 2;;
     elif k=-1 and s=0 and t=1 then return 2;;
     elif k=1 and s=0 and t=0 then return 5;;
     elif k=-1 and s=0 and t=0 then return 5;;
     fi;;
end;;
Q := function(k,s,t) local q;;
     q := function(j,k,s,t)
          if j=1 then return Position(List([1..r(k,s)],
          i->Bt1[3]*Bc[N(k,s)[i]]-t*Bc[N(k,s)[i]]),0*x[1]);;
          else return q(j-1,k,s,t)+Position(List([q(j-1,k,s,t)+1..r(k,s)],
          i->Bt1[3]*Bc[N(k,s)[i]]-t*Bc[N(k,s)[i]]),0*x[1]);;
          fi;;
     end;;
     return List([1..m(k,s,t)],i->N(k,s)[q(i,k,s,t)]);;
end;;
W := [[1,1,-1],[1,1,0],[1,-1,1],[1,-1,0],[1,0,1],[1,0,-1],[1,0,0]];;
S := 2 + 2*Sum(List([1..Length(W)],i->Length(Q(W[i][1],W[i][2],W[i][3]))))
     + Dimension(le)-Dimension(c);;
A := function(i) return List([1..m(W[i][1],W[i][2],W[i][3])],
     t->List([1..m(W[i][1],W[i][2],W[i][3])],
     s->Bc[Q(W[i][1],W[i][2],W[i][3])[s]]*Bc[Q(-W[i][1],-W[i][2],-W[i][3])[t]]));;
end;;
# A(1), A(2), A(3), A(5), A(6) are nonsingular.
# A(7) is singular of order 5 of rank 4; its minor
List([1..4],s->List([1..4],
t->Bc[Q(W[7][1],W[7][2],W[7][3])[s]]*Bc[Q(-W[7][1],-W[7][2],-W[7][3])[t]]));;
# is different from 0.
a := Bc[M(2)[1]]*Bc[M(-2)[1]];;
\end{verbatim}
In conclusion, Condition (2) of Proposition~\ref{p2ar2} is satisfied for
${\goth t}:=\k$\texttt{t} and ${\goth t_1}:=$span(\texttt{Bt1}).

\end{document}